\documentclass[10pt,a4paper]{article}
\usepackage{amsmath}
\usepackage{amssymb}
\usepackage{theorem}
\usepackage{xspace}
\usepackage{amsfonts}
\usepackage{titlesec}
\titleformat{\section}
        {\normalfont\filright\normalsize} 
        {\thesection}
        {0.8em}
        {\bf}[]
        
    \titleformat{\subsection}
        {\normalfont\filright\normalsize} 
        {\thesubsection}
        {0.8em}
        {\bf}[]
\usepackage{epsfig}
\usepackage{hhline}
\usepackage[utf8]{inputenc}    % pour les accents
\theoremstyle{plain}
\newtheorem{thm}{Theorem}[section]
\newtheorem{lem}{Lemma}[section]
\newtheorem{rmk}{Remark}[section]
\newtheorem{prp}{Proposition}[section]
\usepackage[runin]{abstract}
\numberwithin{equation}{section} 
\numberwithin{table}{section} 
\numberwithin{figure}{section}
\def\argmin{\mathop{\mathrm{arg\,min}}}

\def\hh{\hspace*{0.5 cm}}

\def\eb{\textrm{\mathversion{bold}$\mathbf{\beta}$\mathversion{normal}}}  \def\del{\textrm{\mathversion{bold}$\mathbf{\gamma}$\mathversion{normal}}}

\def\ebo{\textrm{\mathversion{bold}$\mathbf{\beta^0}$\mathversion{normal}}}
\def\el{\textrm{\mathversion{bold}$\mathbf{\lambda}$\mathversion{normal}}}
\def\ephi{\textrm{\mathversion{bold}$\mathbf{\phi}$\mathversion{normal}}}
\def\eephi{\textrm{\mathversion{bold}$\Delta$\mathversion{normal}}}
\def\epsi{\textrm{\mathversion{bold}$\mathbf{\psi}$\mathversion{normal}}}

\def\eR{I\!\!R}
\def\eE{I\!\!E}
\def\eP{I\!\!P}
\def\e1{1\!\!1}

\def\ef{\mathbf{\overset{.}{f}}}
\def\eff{\mathbf{\overset{..}{f}}}

\def\eg{\mathbf{{g}}}
\def\egg{\mathbf{\overset{.}{g}}}
\def\eV{\mathbf{{V}}}
\def\es{\mathbf{{S}}}

\def\ev{\mathbf{{v}}}

\def\et{\mathbf{{z}}}
\def\eh{\mathbf{{H}}}
\def\em{\mathbf{{M}}}

\def\ett{\mathbf{\overset{.}{z}}}

\def\XX{\textrm{\mathversion{bold}$\mathbf{X}$\mathversion{normal}}}
\def\xx{\textrm{\mathversion{bold}$\mathbf{x}$\mathversion{normal}}}
\def\ell{\textrm{\mathversion{bold}$\mathbf{\overset{.}{\lambda}}$\mathversion{normal}}}
\newcommand{\R}{\mathbb{R}}
\begin{document}
\title {{\Large
 \bf Empirical likelihood confidence regions for the parameters of a two phases nonlinear model with and without missing response data}}
\author{{Zahraa Salloum } 
\\
\\
{\it  \footnotesize
Université de Lyon, Université Lyon 1, CNRS, UMR 5208, Institut Camille Jordan, Bat. Braconnier,} \\
{\it  \footnotesize
 43, blvd du 11 novembre 1918, F - 69622 Villeurbanne Cedex, France} \\
 {\it  \footnotesize
 Email: salloum@math.univ-lyon1.fr}
}
\maketitle
\rule{\linewidth}{.5pt}
\section*{ \normalsize
\textbf{Abstract}}
\hh In this paper, we use the empirical likelihood method to construct the confidence regions for the
difference between the parameters of a two-phases nonlinear model with random design. We show that the empirical likelihood ratio has an asymptotic chi-squared distribution. The result is a nonparametric version of Wilk's theorem. Empirical likelihood method is also used to construct the confidence regions for the difference between the parameters of a two-phases nonlinear model  with response variables missing at randoms (MAR). In order to construct the confidence regions of the parameter in question, we propose three empirical likelihood statistics : Empirical likelihood based on complete-case data, weighted empirical likelihood and empirical likelihood with imputed values. We prove that all three empirical likelihood ratios have asymptotically chi-squared distributions. The effectiveness of the proposed approaches in aspects of coverage probability and interval length is demonstrated by a Monte-Carlo simulations.
\\ \\
{\it Key-words:} Confidence region; Empirical likelihood; Two-phases problem; Nonlinear regression model; Missing response.
\\
{\it AMS Subject Classification:}  62F03; 62G20; 62J02. \\
%%%%%%%%%%%%%%ùù
%%%%%%%%%%%%%%%%%%%%
\rule{\linewidth}{.5pt} %% Ligne
%%%%%%%%%%%%%%%%ù
%%%%%%%%%%%%%%%%
\section{Introduction}  
Let us consider the following nonlinear model 
\begin{equation}
\label{eq1} 
Y_i=\left\{ 
\begin{array}{ccl}
f(\XX_i;\eb)+\varepsilon_i &  & i=1, \ldots, k, \\
f(\XX_i;\eb_1)+\varepsilon_i &  & i=k+1, \ldots, n,
\end{array}
\right.
\end{equation}
where $\eb$ and $\eb_1$ are $d\times1$ vectors of unknown parameters, $\XX_i$ is a $(p\times 1)$ random vector of regressors with  distribution function H(\textbf{x}), for $\xx \in \Upsilon  $ and $\Upsilon \subseteq \R^p$ a compact set. Let us consider the vector $\textbf{Y}=(Y_1, \ldots, Y_n)$, where, for each observation $i$, $Y_i$ denotes the response variable (which can have missing value) and $\varepsilon_i$ is the error. The continuous random vector sequence $(\XX_i,\varepsilon_i)_{1 \leq i \leq n}$ is independent identically distributed (i.i.d), with the same joint distribution as $(\XX,\varepsilon)$. For all $i$, $\varepsilon_i$ is independent of $\XX_i$.\ \\ 
For the model (\ref{eq1}), let us consider following parameter the difference between the parameters of the two-phases of model $\del=\eb-\eb_1$. \ \\ 
There are two aims in this paper. First, we suppose that in model (\ref{eq1}), function $f$ is non-linear in $\eb$ and the response variable $Y_i$ is observed for each observation $i$. We construct the asymptotic confidence region for $\del$, or we test the null hypothesis
$$H_0 : \del=\del_0,$$
with $\del_0$ a known vector. \\
Second, we construct the confidence regions for $\del$, or we test $H_0$ when some values of $\textbf{Y}$ may be missing and $\XX_i$ is observed completely. That is, we obtain an incomplete sample $\{( \XX_i, Y_i,\delta_i)_{1 \leq i \leq n}\}$ from model (\ref{eq1}), where all the $\XX_i$ are observed, $(\delta_i)_{1 \leq i \leq n}$ being a sequence of random variables, such that $\delta_i=0$ if $Y_i$ is missing and $\delta_i=1$ otherwise. We assume that $Y_i$ is missing at random (MAR). The MAR assumption implies that $\delta_i$ and $Y_i$ are conditionally independent given $\XX_i$. That is, $\eP[\delta_i=1 |\XX_i, Y_i] =\eP[\delta_i=1 |\XX_i]$, for all $1\leq i \leq n$. The MAR assumption is a common condition for statistical analysis with missing data and is reasonable in many practical situations, see Little and Rubin (1987), Qin et al. (2009) and Ciuperca (2011). \ \\

In this kind of problem, we can use the bootstrap approach to construct confidence regions for $\del$, but, one of the inconvenience of the bootstrap is that, it needs some subjective instructions on the shapes and orientations of the confidence regions. In this paper, we will apply the empirical likelihood method for constructing the confidence regions nonparametrically, as an alternative to the bootstrap method. An important characteristic of empirical likelihood is that, it uses only the data to determine the shape and orientation of a confidence regions. This method was introduced by Owen (1988,1990) as a way to extend the ideas of likelihood based inference to certain nonparametric situations.\ \\

Various authors extend empirical likelihood methodology to many statistical situations. To construct the confidence regions for the coefficients in the linear regression model, Chen (1994) proposed a nonparametric method based on empirical likelihood. Qin et al. (2009), Xue (2009) and Ciuperca (2011) considered this same problem but for the models with missing response data. Kolaczyk (1994) shows that empirical likelihood is justified as a method of inference for a class of linear models, and shows in particular how empirical likelihood may be used with generalized linear models. For models with change-points, Kim and Siegmund (1989), Liu et al. (2008) and Ciuperca and Salloum (2013) used the empirical likelihood to detect the change-point in the regression parameters of the linear and nonlinear model. For a epidemic change model, Ning et al. (2012) proposed a method based on the empirical likelihood to detect the epidemic changes of the mean after unknown change points. Always using the empirical likelihood method, Zi et al. (2010) construct the confidence regions for the difference in value between coefficients of two-sample linear regression model with complete data and Wei et al. (2013) for a model with missing response data.\ \\

In this paper, for the model (\ref{eq1}), we use the empirical likelihood method to construct the asymptotic confidence region for $\del=\eb-\eb_1$, firstly, if the response variable $Y_i$ is observed for each $i=1,\ldots,n$, next when the response variable $Y_i$ can be missing. For a model with complete data, we propose an empirical likelihood statistic and we prove that it has a chi-squared asymptotic distribution, which will imply the asymptotic confidence region for $\del$. For a model with missing response data we propose three test statistics : Empirical likelihood based on complete-case data, weighted empirical likelihood and empirical likelihood with imputed values and we show that all three empirical likelihood ratios have asymptotically chi-squared distributions. Then, we generalize the papers of Zi et al. (2010) and Yu et al. (2013) in the nonlinear model case. One of the major  difficulties  for  nonlinear model (beside the linear model approach) is that, for finding the test statistic, the corresponding score functions  depend on the regression parameters, and above all, the analytical form of these derivatives is unknown. On the other hand, in the linear models, many proofs are based on the convexity of the regression function with respect to the parameter regression, then, the extreme value of a convex function is attained on the boundary. These two factors lead to a more difficult theoretical study of the test statistics for nonlinear model.\ \\

The paper is organized as follows. In Section 2, we introduce assumptions, some notations, null and alternative hypothesis. In Section 3, for the model with complete data, we formulate the empirical likelihood ratio and we prove that the empirical likelihood statistic has a chi-squared distribution asymptotically and based to this, we construct the asymptotic confidence region for $\del$. 
%the confidence regions of $\del$ and we give the Wilks theorem for empirical likelihood statistic. 
For a model with missing response data, the asymptotic confidence region of $\del$ corresponding to three proposed empirical likelihood statistics are given in Section 4. The simulation results are presented in Section 5. Proofs of the main results and lemmas are given in Section 6.
\section{Hypothesis, notations, assumptions}
All vectors are column and $\textbf{v}^t$ denotes the transposed of $\textbf{v}$. All vectors and matrices are in bold. Concerning the used norms, for a m-vector $\ev=(v_1, \ldots,v_m)$, let us denote by $\|\ev\|_1= \sum   _{j=1} ^ m |v_j|$ its $L_1$-norm and $\|\ev\|_2=(\sum   _{j=1} ^ m v^2_j)^{1/2}$ its $L_2$-norm. For a matrix $ \textbf{D}=(a_{ij})_{\substack{1\leqslant i \leqslant m_1\\1 \leqslant j \leqslant m_2}}$, we denote by $ \| \textbf{D}\|_1 = \max   _ {j=1,\ldots, m_2} (\sum  _{i=1}^{m_1} |a_{ij}|)$, the subordinate norm to the vector norm $\| .\|_1$. Let $ \overset{\cal L} {\underset{n \rightarrow \infty}{\longrightarrow}}$, $ \overset{\eP} {\underset{n \rightarrow \infty}{\longrightarrow}}$, $ \overset{a.s.} {\underset{n \rightarrow \infty}{\longrightarrow}}$ represent convergence in distribution, in probability and almost sure, respectively, as $n \rightarrow \infty$.   \\
All throughout the paper, C denotes a positive generic constant which may take different values in different formula or even in different parts of the same formula. Moreover, $\textbf{0}_d$ and $\textbf{1}_d$ denote the $d$-vectors with all components zero and 1, respectively. \\

For the model (\ref{eq1}), the regression function $f:\Upsilon \times \Gamma \rightarrow \R$, with  $\Upsilon \subseteq \R^p$ and $\Gamma \subseteq \R^d$, is known up to a $d$-dimensional  parameter $\eb=(\beta_1,\ldots,\beta_d)$. The sets $\Upsilon$ and $\Gamma$ are compact.  \ \\

We now state the assumptions on the errors, on the design and on the regression function. With regard to the random variable $\varepsilon$ we make following assumption :\\
{\bf (A1)} $\eE[\varepsilon_i]=0$, $\eE[\varepsilon_i^2] = \sigma_1^2 < \infty$, for $i=1, \ldots,k$ and $\eE[\varepsilon_j]=0$, $\eE[\varepsilon_j^2] = \sigma_2^2 < \infty$, for $j=k+1, \ldots,n$.\\
The regression function $f:\Upsilon \times \Gamma \rightarrow \R$ and the random vector $\XX$ satisfy the conditions :\\
\textbf{(A2)} for all $\xx \in \Upsilon$ and for $\eb \in \Gamma$, the function $f(\xx,\eb)$ is thrice differentiable in $\eb$ and continuous on $\Upsilon$.\\
In following, for $\xx \in \Upsilon$ and  $\eb \in \Gamma$, we use notation $\ef(\xx,\eb) \equiv \partial f(\xx,\eb)/\partial \eb$, $\eff(\xx, \eb)\equiv \partial^2 f(\xx,\eb)/\partial \eb^2$. \ \\ 
\noindent \textbf{(A3)} $(| \frac{\partial ^2 f(\xx,\eb)}{\partial \beta_s \partial \beta_r }|)_{1 \leq s,r \leq d}$ and $(| \frac{\partial ^3 f(\xx,\eb)}{\partial \beta_s \partial \beta_r \partial \beta_l}|)_{1 \leq s,r,l \leq d}$ are bounded for any $\xx\in \Upsilon$ and $\eb$ in a neighborhood of $\eb^0$. \\ 
\noindent  \textbf{(A4)} $\eE[\| \ef(\XX,\eb) \|_1]< \infty$, 
$\eE[\| \ef(\XX,\eb)\ef^t(\XX,\eb)\|_1]<\infty$ and $\eE[ | \frac{ \partial ^2 f(\XX,\eb) }{\partial \beta_s \partial \beta_r}| ]<\infty$, for all $ 1 \leq s,r \leq d$ and $\eb$ in a neighborhood of $\eb^0$. \\ \\
Assumptions (A3) and (A4) are standard conditions, which are used in nonlinear models, see the book of Seber and Wild (2003) and the paper of Ciuperca and Salloum (2013) for example.\\

We are interested in constructing asymptotic confidence region for $\del=\eb-\eb_1$ and at the same time, in testing the hypothesis 
\begin{equation}
\label{H0}
H_0 : \del=\del_0,
\end{equation}
where $\del_0$ is a $(d \times 1)$ known vector. The alternative hypothesis of (\ref{H0}), is
\begin{equation}
\label{H1}
H_1 : \del \neq \del_0.
\end{equation}
\noindent Under hypothesis $H_0$, let $\ebo$ denote the true value (unknown) of $\eb$, where $\eb$ is the generic value of the regression parameter for the first phase and $\eb_1^0$ the true value of $\eb_1$ for the second phase. If $\del_0^0= \textbf{0}_d$, then $\eb^0_1=\ebo$.\\ 
The change-point location k depends on n, but to simplify notations we use k everywhere. On $k$, we make the classical assumption that $ \lim _{n \rightarrow \infty} \frac{k}{n} \in (0,1)$, see the paper of Ciuperca and Salloum (2013) for example. On the other hand, the change-point location $k$ is supposed fixed for $n$ given. Also to simplify notation, we consider the following sets $I \equiv \{1, ..., k \}$ and $J\equiv \{ k+1, ..., n \}$, which contain the observation subscripts of the two segments for the model (\ref{eq1}). \\
For $i \in I$, let us consider the following $d$-random vectors
$$\eg_i(\eb) \equiv \ef(\XX_i, \eb) [ Y_i - f(\XX_i, \eb)]. $$
For the second phase, we will not consider $Y_j$, but a new variable, noted $Y^*_j$ calculated on the basis of $Y_j$ and taking into account the difference between the parameters $\eb-\del_0$ for the regression function. So, for $ j \in J$, let be $Y^*_j  \equiv Y_j-f(\XX_j, \eb-\del^0)+f(\XX_j, \eb)$ and 
$$\eg_j(\eb) \equiv \ef(\XX_j, \eb) [Y^*_j - f(\XX_j, \eb)].$$ We remark that, under hypothesis $H_0$, we have $\eg_i(\ebo)=\ef(\XX_i,\ebo) \varepsilon_i$ for $ i \in I$ and $\eg_j(\ebo)=\ef(\XX_j,\ebo) \varepsilon_j$ for $ j \in J$. For all $i=1, \ldots, n$, we have $\eE[\eg_i(\ebo)]=\textbf{0}_d$. \\
We consider also the $d \times d$ matrix 
$$\eV \equiv \eE[ \ef(\XX_i,\ebo) \ef^t(\XX_i,\ebo)].$$
\hh In order to introduce the maximum empirical likelihood method in the following section, let $y_1, \ldots, y_k,$ \\$ y_{k+1}, \ldots, y_n$ be observations for the random variables $Y_1, \ldots, Y_k, Y_{k+1}, \ldots, Y_n$. Corresponding to the sets $I$ and $J$, let be the probability vectors $(p_1,\ldots,p_k)$ and $(q_{k+1},\ldots,q_{n})$. These vectors contain the probability to observe the value $y_i$ (respectively $y_j$) for the dependent variable $Y_i$ (respectively $Y_j$) : $p_i\equiv \eP [Y_i=y_i]$, for $i=1, \ldots, k$ and $q_j \equiv \eP[Y_j=y_j]$, for $j=k+1, \ldots,n$. Obviously, these probabilities satisfy the relations $\sum _{i \in I}p_i = 1$ and $\sum _{j \in J}q_j = 1$. 
%We consider also for $i\in I, \Var (\varepsilon_i\ef(\XX_i,\ebo))= \sigma_1^2\eV$, and for $j\in J, \Var (\varepsilon_j\ef(\XX_j,\ebo))= \sigma_2^2\eV$. \\ \\
\section{Model with complete data}
\label{sec3}
In this section, we suppose that, for the nonlinear model given by (\ref{eq1}), the response variable $Y_i$ is observed for each $i=1,\ldots,n$. We will construct the empirical likelihood ratio statistic and show that this statistic has a $\chi^2$ asymptotic distribution, which allows us to construct the asymptotic confidence region for $\del$.
\subsection{Test statistic}
\hspace*{0.3cm} In this subsection, we formulate the empirical likelihood ratio statistic which will be used to construct the asymptotic confidence region for $\del=\eb-\eb_1$, or for testing hypothesis $H_0$, given by (\ref{H0}), against the alternative $H_1$, given by (\ref{H1}). \ \\

Under  hypothesis $H_0$, we have $\del_0=\ebo-\eb_1^0$. Remark that if $\del_0$, $\ebo$ are given then $\eb^0_1$ (the true value of $\eb_1$ under $H_0$) is known. In order to study the empirical likelihood statistic, since $\ebo$ is unknown we will use the notation $\eb$. Then, the profile empirical likelihood for $\del$, evaluated at $\del_0$ under $H_0$ is defined as
\begin{eqnarray*}
{\cal R}_{nk}(\del_0,\eb)&\equiv  &\sup   _{(p_1,\ldots,p_k)}
  \sup   _{(q_{k+1},\ldots,q_n)}  \bigg\{  \prod    _{i \in I}p_i   \prod    _{j \in J}q_j ;   \sum    _{i\in I} p_i = 1,  \sum    _{j\in J} q_j =1, 
  \nonumber \\ &&
 \sum    _{i \in I}p_i \eg_i( \eb )= \sum    _{j \in J}
q_j \eg_j( \eb)=\textbf{0}_d \bigg \}.
\end{eqnarray*} 
%%%%%%%%%%%%%%%%%%%%%%%%%%%%%%%%%%%%%%%%%%%%%%%%%%%%%%%%
Without constraints $ \sum  _{i\in I} p_i \eg_i(\eb)=\textbf{0}_d$, the maximum of $\prod   _{i\in I} p_i$, $\prod   _{j \in J}q_j$ are  attained for $p_i=k^{-1}$, $q_j=(n-k)^{-1}$, respectively. Then, the profile empirical likelihood ratio for $\del$, evaluated at $\del_0$ under $H_0$ has the form 
\begin{eqnarray*}
\begin{array}{ccl}
\label{R'0}
{\cal R}'_{nk}(\del_0,\eb) \equiv  \sup  \limits_{(p_1,\cdots,p_k)}
  \sup  \limits  _{(q_{k+1},\cdots,q_n)} \bigg \{  \prod  \limits  _{i \in
I}k p_i   \prod  \limits_{j \in J}(n-k) q_j ;   \sum \limits   _{i\in I} p_i = 1,  \sum  \limits  _{j\in J} q_j =1,\\   \sum  \limits  _{i \in I}p_i \eg_i( \eb )= \sum   \limits  _{j \in J}q_j \eg_j( \eb)=\textbf{0}_d  \bigg\}.
\end{array}
\end{eqnarray*}
%%%%%%%%%%%%%%%%%%%%%%%%%%%%%%%%%%%%%%%%%%%%%%%%%%%%%
Under hypothesis $H_1$, we can consider estimators for the parameters $\eb$ and $\eb_1$. Then, let us consider the least square estimators $\hat{\eb}$ and $\hat{\eb}_1$ of $\eb$ and $\eb_1$, on the observations corresponding to the sets $I$ and $J$
\begin{equation}
\label{eqb1b2}
 \hat{\eb}= \argmin \limits_{\eb} \sum _{i \in I}(Y_i-f(\XX_i,\eb))^2, \,\,\,\,
\hat {\eb}_1=\argmin \limits_{\eb_1} \sum _{j \in J}(Y_j-f(\XX_j,\eb_1))^2.
\end{equation}
%%%%%%%%%%%%%%%%%%%%%%%%%%%%%%%%%%%%%%%%%%%%%%%%%%%%%
Then, the profile empirical likelihood for $\del$ under hypothesis $H_1$ is
\begin{eqnarray*}
{\cal R}_{nk}(\hat \del,\hat\eb, \hat\eb_1)&\equiv  &\sup   _{(p_1,\ldots,p_k)}
  \sup   _{(q_{k+1},\ldots,q_n)}  \bigg\{  \prod    _{i \in I}p_i   \prod    _{j \in J}q_j ;   \sum    _{i\in I} p_i = 1,  \sum    _{j\in J} q_j =1, 
  \nonumber \\ &&
 \sum    _{i \in I}p_i \eg_i( \hat \eb )= \sum    _{j \in J}
q_j \eg_j( \hat \eb_1)=\textbf{0}_d \bigg \},
\end{eqnarray*} 
where $\hat{\del}=\hat{\eb}-\hat{\eb}_1$, with $\hat{\eb}$, $\hat {\eb}_1$ given by equation (\ref{eqb1b2}). \\
%%%%%%%%%%%%%%%%%%%%%%%%%%%%%%%%%%%%%%%%%%%%%%%%%%%%%
The profile empirical likelihood ratio for $\del$, under hypothesis $H_1$ has the form 
\begin{eqnarray*}
\begin{array}{ccl} 
{\cal R}'_{nk}(\hat \del,\hat\eb, \hat\eb_1)\equiv  \sup \limits  _{(p_1,\cdots,p_k)}\sup \limits _{(q_{k+1},\cdots,q_n)}
\bigg \{  \prod  \limits  _{i \in I}k p_i \prod \limits  _{j \in J}(n-k) q_j
 ;   \sum  \limits _{i\in I} p_i =1, \sum \limits  _{j\in J} q_j =1, \\   \sum \limits_{i \in I}p_i \eg_i( \hat \eb )=\textbf{0}_d,\sum \limits  _{j\in J}q_j \eg_j( \hat \eb_1)=\textbf{0}_d \bigg \}.
 \end{array}
\end{eqnarray*}
%%%%%%%%%%%%%%%%%%%%%%%%%%%%%%%%%%%%%%%%%%%%%%%%%%%%%%%
Thus, using an idea similar to the maximum likelihood test for testing $H_0$ against $H_1$, we consider the profile empirical likelihood ratio
\begin{equation}
\label{EL}
\frac{{\cal R}_{nk}(\del_0,\eb)}{{\cal R}_{nk}(\hat{\del}, \hat{\eb}, \hat{\eb}_1)}=\frac{{\cal R}'_{nk}(\del_0,\eb)}{{\cal R}'_{nk}(\hat{\del}, \hat{\eb}, \hat{\eb}_1)}.
\end{equation}
By Theorem 1 of Ciuperca and Salloum (2013), under hypothesis $H_1$ the profile empirical likelihood ratio for $\hat \eb$, $ \hat \eb_1$ has a $\chi ^2$ asymptotic distribution with $2d$ degrees of freedom. Then, ${\cal R}'_{nk}(\hat \del,\hat\eb, \hat\eb_1)$ is not asymptotically depend on the parameters $\hat \eb$ and $ \hat \eb_1$. Then, the corresponding empirical log-likelihood function denoted by ${\cal \tilde Z}_{nk}(\del_0,\eb)$ can be written as
 \begin{eqnarray}
 \label{z1}
{\cal \tilde Z}_{nk}(\del_0,\eb)
%&=& -2 \log \Big[{\cal R}_{nk}(\del_0,\eb)/{\cal R}_{nk}(\hat{\del}, \hat{\eb}) \Big] 
&=& \sup   _{(p_1,\ldots,p_k)}\sup  _{(q_{k+1},\ldots,q_n)}
\bigg \{  \sum    _{i \in I}\log (k p_i )+\sum   _{j \in J} \log ((n-k) q_j)
 ;   \sum   _{i\in I} p_i =1, 
 \nonumber \\ && 
\qquad  \sum   _{j\in J}  q_j =1, \sum_{i \in I}p_i \eg_i( \eb)=\textbf{0}_d,\sum   _{j\in J}q_j \eg_j( \eb)=\textbf{0}_d \bigg \}.
\end{eqnarray}
In order to study ${\cal \tilde Z}_{nk}(\del_0,\eb)$ of (\ref{z1}), consider first the observations on the first phase, for $i \in I$. Using the Lagrange multiplier method, the following random process 
$$\sum  _{i \in I} \log p_i + \eta (\sum _{i \in I} p_i -1)-k \el_1^t \sum  _{i \in I}  p_i \eg_i(\eb),$$
where $\el_1$ and $\eta$ are the Lagrange multipliers, $\el_1 \in \R^d$ and $\eta \in \R$, must be maximized with respect to $p_1, \ldots,p_k$, $\eta$ and $\el_1$.
%Taking derivative with respect to $p_i$ of this process equal to zero, we obtain
Using the fact that the derivative of the process given above with respect to $p_i$ is equal to zero, we obtain
\begin{equation}
\label{new}
p_i= \frac{1}{k \el_1^t \eg_i(\eb)- \eta}.
\end{equation}
By relation (\ref{new}), we can obtain easily that $\eta= k$. Hence, the probability $p_i$ becomes
\begin{equation}
\label{eqq1}
p_i=\frac{1}{k(1+\el^t_1 \eg_i(\eb))}.
\end{equation}
In the same way, we obtain for the observations on the second phase of (\ref{z1}) for $j \in J$,
\begin{equation}
\label{eqq2}
q_j=\frac{1}{n-k-n\el^t_2 \eg_j(\eb)},
\end{equation}
where $\el_2\in \R^d $ is the Lagrange multiplier. \\
Using equations (\ref{eqq1}) and (\ref{eqq2}), to study statistic of (\ref{z1}) amounts to maximizing, with respect to $\el_1$ and $\el_2$, the following random process
\begin{equation}
\label{z2}
{\cal \tilde Z}_{nk}(\del_0,\el_1,\el_2,\eb)\equiv 2 \Big [  \sum \limits _{i \in I} \log(1+\frac{n}{k} \el_1^t\eg_i(\eb))+\sum \limits _{j \in J} \log(1-\frac{n}{n-k} \el_2^t\eg_j(\eb))\Big].
\end{equation}
In order to have single Lagrange multiplier, denoted by $\el$, we restrict the study to a particular case, when $\el_1$ and $\el_2$ satisfy the constraint $\eV_{1n}(\eb) \el_1 =\eV_{2n}(\eb) \el_2 $, with 
$$\eV_{1n}(\eb) \equiv k^{-1} \sum  _{i \in I} \egg_i(\eb), \quad \eV_{2n}(\eb) \equiv (n-k)^{-1} \sum   _{j \in J} \egg_j(\eb).$$
In the case of the true parameter $\ebo$, this two last matrices are denoted $\eV^0_{1n} \equiv \eV_{1n}(\ebo)$ and $\eV^0_{2n} \equiv \eV_{2n}(\ebo)$. \\
Considering constraint  $\eV_{1n}(\eb) \el_1 =\eV_{2n}(\eb) \el_2 $, statistic (\ref{z2}) becomes
\begin{equation}
\label{z3}
{\cal \tilde{Z}}_{nk}(\del_0,\el,\eb)\equiv 2 \Big[  \sum \limits _{i \in I} \log(1+\frac{n}{k} \el^t\eg_i(\eb))+\sum \limits _{j \in J} \log(1-\frac{n}{n-k} \el^t \eV_{1n}(\eb)\eV_{2n}^{-1}(\eb) \eg_j(\eb)) \Big].
\end{equation}
We will study the maximum, with respect to $\eb$ and $\el$, of empirical log-likelihood test statistic ${\cal \tilde{Z}}_{nk}(\del_0,\el,\eb)$. Then, we calculate the score functions of test statistic (\ref{z3})
\begin{eqnarray} 
\label{phi1}
\tilde {\ephi}_{1n}(\del_0,\el,\eb) & \equiv & \frac{\partial {\cal \tilde{Z}}_{nk}(\el,\del_0,\eb)}{2 \partial \el} 
 \nonumber \\ 
 &=& \sum  _{i \in I}
\frac{\eg_i(\eb)}{\frac{k}{n}+\el^t \eg_i(\eb)}-\sum  _{j\in J}
\frac{\eV_{1n}(\eb)\eV_{2n}^{-1}(\eb)\eg_j(\eb)}{\frac{n-k}{n}-\el^t\eV_{1n}(\eb)\eV_{2n}^{-1}(\eb)\eg_j(\eb)}.
\nonumber \\
\end{eqnarray}
\begin{eqnarray}
\label{phi2}
\tilde {\ephi}_{2n}(\del_0,\el,\eb) & \equiv & \frac{\partial {\cal \tilde{Z}}_{nk}(\el,\del_0,\eb)}{2 \partial \eb} 
 \nonumber \\ 
 &=& \sum  _{i\in I} \frac{
\egg_i(\eb) \el^t}{\frac{k}{n}+
\el^t(\eb)\eg_i(\eb)}  
-\sum   _{j\in
J}  \frac{\partial (\eV_{1n}(\eb) \eV_{2n}^{-1}(\eb)\eg_j(\eb))/ \partial \eb }
{\frac{n-k}{n}-\el^t \eV_{1n}(\eb) \eV_{2n}^{-1}(\eb)\eg_j(\eb)}\el^t. \nonumber \\
\end{eqnarray}
Thus, solving the system $\tilde {\ephi}_{1n}(\del_0,\el,\eb) =\textbf{0}_d $ and $\tilde {\ephi}_{2n}(\del_0,\el,\eb) =\textbf{0}_d$, the obtained solutions $\tilde \el_n(k)$ and $\tilde \eb_n (k)$ are the maximizers of the statistic (\ref{z3}).\ \\ \\
We emphasize that, compared with a linear model, in our case, matrix $\eV_{1n}(\eb)$, $\eV_{2n}(\eb)$ and derivative $\egg(\eb)$ depend on $\eb$. These, besides the nonlinearity of $\eg(\eb)$ involve difficulties in the study of the statistic ${\cal \tilde{Z}}_{nk}(\del_0,\el,\eb)$ and of the solutions $\tilde \el_n(k)$, $\tilde \eb_n(k)$. 
\begin{rmk} 
\label{rem1} 
To acquire the symmetric form of the statistic ${\cal \tilde{Z}}_{nk}(\del_0,\el,\eb)$ given by (\ref{z3}), which makes the arguments more concise, we consider the following notations, for $ i \in I$ and $j \in J$
\begin{eqnarray*}
\et_i(\eb) \equiv  \em_n^{\frac{1}{2}}(\eb)\eV_{1n}^{-1}(\eb)\eg_i(\eb),  \qquad \et_j(\eb) \equiv  \em_n^{\frac{1}{2}}(\eb)\eV_{2n}^{-1}(\eb)\eg_j(\eb),
 \end{eqnarray*}
 where
 \begin{eqnarray*}
\em_n(\eb) \equiv \frac{k(n-k)}{n^2} \eV_{1n}(\eb) \eh_n(\eb) \eV_{2n}(\eb)
 \end{eqnarray*}
 and
\begin{eqnarray*}
\eh_n(\eb)\equiv  \Big[\frac{k}{n} \sigma^2_2 \eV_{1n}(\eb)+\frac{n-k}{n} \sigma^2_1 \eV_{2n}(\eb)\Big]^{-1}. 
\end{eqnarray*}
Taking into account the above notations, we consider instead of ${\cal \tilde Z}_{nk}(\del_0,\el,\eb)$ given by (\ref{z3}), the following test statistic
\begin{equation}
\label{zfin}
{\cal Z}_{nk}(\del_0,\el,\eb) \equiv 2 \Big[  \sum \limits _{i \in I} \log(1+\frac{n}{k} \el^t\et_i(\eb))+\sum \limits _{j \in J} \log(1-\frac{n}{n-k} \el^t \et_j(\eb) \Big].
\end{equation}
Based on (\ref{zfin}), we can derive the following score equations to get the estimators $(\hat \el_n(k),\hat \eb_n(k))$ of $(\el,\eb)$ 
\begin{equation}
\label{scoreZ}
\left\{
\begin{array}{ccl} 
{\ephi}_{1n}(\del_0,\el,\eb) \equiv \displaystyle \frac{\partial {\cal Z}_{nk}(\del_0,\el,\eb)}{2 \partial \el} = \sum \limits _{i \in I}
\frac{\et_i(\eb)}{\frac{k}{n}+\el^t \et_i(\eb)}-\sum \limits _{j\in J}
\frac{\et_j(\eb)}{\frac{n-k}{n}-\el^t\et_j(\eb)}, \\ \\ 
{\ephi}_{2n}(\del_0,\el,\eb) \equiv \displaystyle \frac{\partial {\cal Z}_{nk}(\del_0,\el,\eb)}{2 \partial \eb}   = \sum \limits _{i \in I}
\frac{\ett_i(\eb)}{\frac{k}{n}+\el^t \et_i(\eb)}-\sum \limits _{j\in J}
\frac{\ett_j(\eb)}{\frac{n-k}{n}-\el^t\et_j(\eb)},
\end{array}
\right.
 \end{equation}
where $\ett_i(\eb)$ and $\ett_j(\eb)$ are the derivative with respect to $\eb$ of $\et_i(\eb)$ and $\et_j(\eb)$ respectively.
\end{rmk}
\hh In applications, the error variance $\sigma^2_1$ can be estimated by $k^{-1} \sum   _{i \in I} [Y_i - f(\XX_i,\hat{\eb})]^2$, and $\sigma^2_2$ can be estimated by $(n-k)^{-1} \sum   _{j \in J} [Y_j  - f(\XX_j,\hat{\eb}_1)]^2$, with $\hat{\eb}$ and $\hat{\eb}_1$ given by relation (\ref{eqb1b2}). 
%%%%%%%%%%%%%%%%%%%%%%%%%%%%%%%%%%%%%%%%%%%%%%%%%%%%%%%%%%%%%%%%%%%%%%%%
%%%%%%%%%%%%%%%%%%%%%%%%%%%%%%%%%%%%%%%%%%%%%%%%%%%%%%%%%%%%%%%%%%%%%%%%%%
%%%%%%%%%%%%%%%%%%%%%%%%%%%%%%%%%%%%%%%%%%%%%%%%%%%%%%%%%%%%%%%%%%%%%%%%%%%%
\subsection{Asymptotic behaviour of the statistic ${\cal Z}_{nk}$}
\hh In this subsection, we will study the asymptotic behaviour of the statistic ${\cal Z}_{nk}$ given by equation (\ref{zfin}) under null hypothesis $H_0$, given by (\ref{H0}). We show also that $\hat \el_n(k)$ and $\hat \eb_n(k)$, the solutions of the score equations $\ephi_{1n}(\del_0,\el,\eb)=\textbf{0}_d$ and $\ephi_{2n}(\del_0,\el,\eb)=\textbf{0}_d$ given by relation (\ref{scoreZ}) have suitable properties. \ \\ \ \\
\hh In addition to assumptions (A1)-(A4), we require the following assumptions : \ \\ \ \\
\noindent  \textbf{(A5)} The matrices $\eV_{1n}(\eb)$ and $\eV_{2n}(\eb)$ are non singular for any $\XX_i \in \Upsilon$ and all $\eb$ in a neighborhood of $\eb^0$, moreover their determinants are bounded for sufficiently large n. \ \\
\noindent  \textbf{(A6)} $ \,\, \sup  \limits _ {\eb \in \Gamma} \eE [ \| \ef (\XX,\eb)] \|_2] ^{2s} < \infty$ and for all $ 1 \leq u,v \leq d$,  $\,\, \sup  \limits_ {\eb \in \Gamma} \eE [ | \frac{\partial ^2 f(\XX,\eb) }{\partial \beta_u \partial \beta_v}|^s] $ $< \infty$, for some $s>2$.\\ 

Assumption (A5) assures that the matrices $\eV_{1n}(\eb)$ and $\eV_{2n}(\eb)$ are uniformly nonsingular and bounded for sufficiently large n. Assumption (A6) is a necessary moment condition for statistical inference and it is also employed in the paper of Boldea and Hall (2013). \ \\

By the next proposition, we show that $\hat \el_n(k)$ and $\hat \eb_n(k)$, the solutions of the score equations $\ephi_{1n}(\del_0,\el,\eb)$ \\$=\textbf{0}_d$ and $\ephi_{2n}(\del_0,\el,\eb)=\textbf{0}_d$ given by relation (\ref{scoreZ}), have suitable properties. More precisely, we show that  $\|\hat \el_n(k)\|_2 \rightarrow 0$, as $n\rightarrow \infty$ and that $\hat \eb_n(k)$ is a consistent estimator of $\eb^0$, under hypothesis $H_0$. The proof is in Appendix, Section 6.
\begin{prp}
\label{prop1}
Under null hypothesis $H_0$, if the assumptions (A1)-(A6) are satisfied, then we have $\hat{\el}_n(k)=O_{\eP}(n^{-1/2})$ and $\hat \eb_n(k)- \eb_0=o_{\eP}(1)$, where $(\hat \el_n(k),\hat \eb_n(k))$ is the solution of the system (\ref{scoreZ}).
\end{prp}
%%%%%%%%%%%%%%%%%%%%%%%%%%%%%%%%%%%%%%%%%%%%%%%%%%%
%%%%%%%%%%%%%%%%%%%%%%%%%%%%%%%%%%%%%%%%%%%%%%%%%%%
The following result is a generalization of the nonparametric version of Wilk's theorem for the empirical likelihood ratio defined by (\ref{zfin}). The proof is in Appendix, Section 6.
\begin{thm}
\label{theo1}
Suppose that assumptions (A1)-(A6) hold. Under null hypothesis $H_0$, the statistic \\
${\cal Z}_{nk}(\del_0,\el,\eb)$ given by (\ref{zfin}) converge, as $n \rightarrow \infty$, to a chi-squared distribution with d degrees of freedom, where $d$ is the dimension of $\eb$.  
\end{thm}
\hh From Theorem \ref{theo1}, for fixed size $\alpha \in (0,1)$, we can construct asymptotic confidence region for $\del$ as follow: $$CR_{\alpha}=\{ \del \in \R^d :{\cal Z}_{nk}(\del_0,\el,\eb)< c_{1-\alpha;d} \},$$
where $c_{1-\alpha;d}$ is the $(1-\alpha)$ quantile of the chi-squared distribution with $d$ degrees of freedom. 
%/************************************************************************
%/************************************************************************
%/************************************************************************
%/************************ Section ************************************
%/************************************************************************
%/************************************************************************
%/************************************************************************
\section{Model with missing response data}
\label{sec4} In this section, for model (\ref{eq1}), we suppose that all the $\XX_i$'s are observed, in exchange the response variable $Y_i$ can be missing. Let the sequence of random variables $(\delta_i)_{1 \leq i \leq n}$ be defined by $\delta_i=0$ if $Y_i$ is missing and $\delta_i=1$ if $Y_i$ is observed. We suppose that $Y_i$ is missing at random (MAR) i.e $\eP[\delta_i=1 |\XX_i, Y_i] =\eP[\delta_i=1 |\XX_i]$, for all $1\leq i \leq n$. \\
We consider the selective probabilities functions defined as $ \pi_1(\textbf{x} _i)= \eP[\delta_i=1 | \XX_i= \textbf{x} _i]$, for $i \in I$ and $ \pi_2(\textbf{x} _j)= \eP[\delta_j=1 | \XX_j= \textbf{x} _j]$, for $j \in J$. We suppose that $ \pi_1(\textbf{x} _i)>0$ and $\pi_2(\textbf{x} _j)>0$, which are a common suppositions in the literature, see for example the papers of Sun et al. (2009), Xue (2009) and Ciuperca (2011). Moreover, $\lim \limits _{n \rightarrow \infty} \sum  _{i \in I} \delta_i = \lim \limits  _{n \rightarrow \infty} \sum  _{j \in J} \delta_j = \infty$, which means that, we have a significant number of values for Y non missing. \\
%The literature on statistical analysis of data with missing values has flourished since the early 1970, spurred by advances in computer-technology that made previously laborious numerical calculations a simple matter. In practice, however, response variables are usually missing due to various reasons such as unwillingness of some sampled units to supply the desired information, loss of information caused by uncontrollable factors, failure on the part of investigators to gather correct information and so forth. Actually, missingness of responses is very common in opinion polls, market research surveys and many scientific experiments. When some responses are missing, the existing methods in the literature are not applicable any more. \\ 

A nonlinear model based on missing at random (MAR), has been considered by various authors. Muller (2009) constructed a efficient estimator for expectation $\eE[h(\XX,\textbf{Y})]$ using a efficient estimator of parameters, with $h$ is a known square integrable function. The mean response $\eE(\textbf{Y})$ is a special case. Ciuperca (2013) constructed the empirical likelihood ratios using complete-case and imputed values. The basic idea in imputation is to "fill in" missing $Y$ values with "appropriate" values to create a completed data set, thereby allowing standard methods to be applied. However, the imputed data are not i.i.d. because a plug-in estimator is used. \ \\ 

For the rest of this section, an empirical likelihood method is used to study model (\ref{eq1}) under missing response data. We are interested to construct the asymptotic confidence region for $\del=\eb-\eb_1$, based on the data $( \XX_i, Y_i,\delta_i)_{1 \leq i \leq n}$, or testing the null hypothesis  
$$H_0:\del= \del_0.$$
We recall that, under hypothesis $H_0$, $\ebo$ denote the true value (unknown) of $\eb$, where $\eb$ is the generic value of the regression parameter for the first phase.
%To construct the confidence regions of $\del$, we propose three empirical likelihood statistics : Empirical likelihood based on complete-case data, weighted empirical likelihood and empirical likelihood with imputed values.  
%%%%%%%%%%%%%%%%%%%%%%%%%%%%%%%%%%%%%%%%%%%%%%%%
%%%%%%%%%%%%%%%%%%%%%%%%%%%%%%%%%%%%%%%%%%%%%%%%%%%%%%%
\subsection{Test statistics}
To construct the asymptotic confidence region of $\del$, we propose three empirical likelihood statistics : empirical likelihood based on complete-case data, weighted empirical likelihood and empirical likelihood with imputed values.  
\subsubsection{Empirical likelihood based on complete-case data} 
\hspace*{0.3cm} Firstly, we give  the empirical likelihood based on complete-case data, i.e, excluding missing data. In the regression context, this
usually means complete-case analysis : excluding all units for which the outcome
or any of the inputs are missing.\\ \\
Two problems arise with complete-case analysis. First, if the units with missing values differ systematically from the completely observed cases, this could bias the complete-case analysis. Second, if many variables are included in a model, there may be very few complete cases, so that most of the data would be discarded for the sake of a simple analysis. 
\ \\

Then, for observations of the response variable of the second phase, we will consider new variable, noted $\tilde Y_j$ defined similarly to complete-case for non missing values. More precisely, for model (\ref{eq1}), for $j \in J$, let be if $Y_j$ non missing, $\tilde Y_j \equiv Y_j-f(\XX_j, \eb-\del^0)+f(\XX_j, \eb)$ and $\tilde Y_j$ equal to any finite value if $Y_j$ missing. Then, we define the two following $d-$random vectors
\begin{equation}
\left \{
\begin{array}{ccl}
\eg_{i,C}(\eb) \equiv \delta_i \ef_i(\eb) (Y_i-f_i(\eb)), &  &  i \in I, \\
\eg_{j,C}(\eb) \equiv \delta_j \ef_j(\eb) (\tilde Y_j-f_j(\eb)), &  &  j \in J.
\end{array}
\right.
\end{equation}
Let us also consider the following matrices
\begin{equation*}
\label{VC}
\eV_{1n,C}(\eb) \equiv \frac{1}{k} \sum  _{i \in I} \egg_{i,C}(\eb), \qquad \eV_{2n,C}(\eb)   \equiv \frac{1}{n-k} \sum   _{j \in J} \egg_{j,C}(\eb).
\end{equation*}
Similarly as in Section 2 and by the same argument given in Remark \ref{rem1}, in order to construct the asymptotic confidence region for $\del$ or to test hypothesis $H_0$ in the complete-case data method, the corresponding empirical likelihood ratio statistic, is
\begin{equation}
\label{zcfin}
{\cal Z}_{nk,C}(\del_0,\el_C,\eb) \equiv 2 \Big[  \sum \limits _{i \in I} \log(1+\frac{n}{k} \el_C^t\et_{i,C}(\eb))+\sum \limits _{j \in J} \log(1-\frac{n}{n-k} \el_C^t \et_{j,C}(\eb) \Big],
\end{equation}
with
\begin{eqnarray*}
\et_{i,C}(\eb)  \equiv  \em_{n,C}^{\frac{1}{2}}(\eb) \eV_{1n,C}^{-1}(\eb)\eg_{i,C}(\eb), \qquad \et_{j,C}(\eb)  \equiv  \em_{n,C}^{\frac{1}{2}} (\eb) \eV_{2n,C}^{-1}(\eb)\eg_{j,C}(\eb),
\end{eqnarray*}
\begin{eqnarray*}
\em_{n,C}(\eb) \equiv \frac{k(n-k)}{n^2} \eV_{1n,C}(\eb) \eh_{n,C}(\eb) \eV_{2n,C}(\eb)
\end{eqnarray*}
and
\begin{eqnarray*}
\eh_{n,C} (\eb) & \equiv & \Big [ \frac{ n-k}{nk} \sigma^2_1 \eV_{2n,C}(\eb) \eV_{1n,C}^{-1}(\eb)  \sum \limits _{i \in I}  \delta_i \pi_1(\textbf{X} _i) \big[\eff_{i}(\eb)(Y_i-f_i(\eb))-\ef_i(\eb)\ef_i^t(\eb) \big]
\\  
  && + \frac{ k}{n(n-k)} \sigma^2_2  \sum \limits _{j \in J}  \delta_j \pi_2(\textbf{X} _j) \big[\eff_j(\eb) (\tilde Y_j-f_j(\eb))-\ef_j(\eb)\ef_j^t(\eb)\big] \eV_{2n,C}^{-1}(\eb) \eV_{1n,C}(\eb) \Big ]^{-1}.
\end{eqnarray*}
We recall that, $\sigma_1^2$ and $\sigma_2^2$ denote the variance of $\varepsilon_i$ and $\varepsilon_j$, respectively. For  model with missing response, under assumptions $\lim \limits_{n \rightarrow \infty} \sum  _{i \in I} \delta_i \rightarrow \infty$ and $\lim \limits_{n \rightarrow \infty} \sum _{j \in J} \delta_j \rightarrow \infty$, they can be estimated respectively by $\tilde \sigma^2_1$ and $\tilde \sigma^2_2$, defined by
\begin{equation}
\label{sigtild}
\tilde{\sigma}_1^2= \frac{\sum \limits  _{i \in I} \delta_i[Y_i-f(\XX_i,\hat{\eb})]^2}{\sum \limits  _{i \in I} \delta_i}, \qquad \tilde{\sigma}_2^2= \frac{\sum \limits  _{j \in J} \delta_j [Y_j-f(\XX_j,\hat{\eb}_1)]^2}{\sum \limits  _{j \in J} \delta_j}, 
\end{equation} 
with
\begin{equation}
\label{eqbb}
\hat{\eb}= \argmin \limits_{\eb} \sum _{i \in I}\delta_i (Y_i-f(\XX_i,\eb))^2, \,\,\,\,
\hat{\eb}_1= \argmin \limits_{\eb_1} \sum _{j \in J} \delta_j (Y_j-f(\XX_j,\eb_1))^2.
\end{equation}
The score functions of test statistic (\ref{zcfin}) are
\begin{equation}
\label{scoreZc}
\left\{ 
\begin{array}{ccl}
\ephi_{1n,C}(\del_0,\el_C,\eb)  = \displaystyle \sum \limits _{i \in I}
\frac{\et_{i,C}(\eb)}{\frac{k}{n}+\el^t_C \et_{i,C}(\eb)}-\sum \limits _{j\in J}
\frac{\et_{j,C}(\eb)}{\frac{n-k}{n}-\el^t_C\et_{j,C}(\eb)}, \\ \\
\ephi_{2n,C}(\del_0,\el_C,\eb)   = \displaystyle \sum \limits  _{i \in I}
\frac{\ett_{i,C}(\eb)}{\frac{k}{n}+\el^t_C \et_{i,C}(\eb)}-\sum \limits _{j\in J}
\frac{\ett_{j,C}(\eb)}{\frac{n-k}{n}-\el^t_C\et_{j,C}(\eb)},
\end{array}
\right.
\end{equation}
where $\ett_{i,C}(\eb)$ and $\ett_{j,C}(\eb)$ are the derivative with respect to $\eb$ of $\et_{i,C}(\eb)$ and $\et_{j,C}(\eb)$ respectively.\\
Then, solving the system ${\ephi}_{1n,C}(\del_0,\el_C,\eb)=\textbf{0}_d$ and ${\ephi}_{2n,C}(\del_0,\el_C,\eb)=\textbf{0}_d$ given by (\ref{scoreZc}), we obtain $\hat \el_{n,C}(k)$ and $\hat \eb_{n,C}(k)$ the maximizers of the statistic (\ref{zcfin}). \ \\ 

So far, the selective probabilities $\pi_1(\XX_i)$ and $\pi_2(\XX_j)$ were considered as known. If they are unknown, we can consider the nonlinear estimators $\hat \pi_1(\XX_i)$ and $\hat \pi_2(\XX_j)$ for $\pi_1(\XX_i)$ and $\pi_2(\XX_j)$, respectively, given by
\begin{equation}
\label{eqq21}
\begin{array}{ccl}
 \hat \pi_1(\XX_i)=  \displaystyle \frac{\sum \limits   _{l \in I} \delta_l K_{1} ((\XX_l-\XX_i )/h_{1n})}{\max \{1,\sum \limits  _{l \in I} K_{1} ((\XX_l-\XX_i )/h_{1n})  \}}, &  & i \in I,  \\ \\
\hat \pi_2(\XX_j)=  \displaystyle \frac{\sum \limits  _{l \in J} \delta_l K_{2} ((\XX_l-\XX_j )/h_{2n})}{\max \{1,\sum \limits  _{l \in J} K_{2} ((\XX_l-\XX_j )/h_{2n})  \}},  &  & j \in J.
\end{array}
\end{equation}
Here, $h_{1n}$ and $h_{2n}$ are a positive sequences tending towards 0 as $n \rightarrow \infty$. $K_{1}$, $K_{2}$ are kernel functions defined in $\eR^d$. \\
The bandwidths $h_{1n}$ and $h_{2n}$ satisfies the following : \\ \\
\noindent \textbf{(A7)} $k h_{1n} ^{4 \max \{2,d-1 \} } \rightarrow 0$ and $(n -k)h_{2n} ^{4 \max \{2,d-1 \} } \rightarrow 0$, as $n \rightarrow \infty$.\\  \\
The kernel functions $K_{1}$ and $K_{2}$ satisfy the classical condition : \\  \\
\noindent \textbf{(A8)} There exist positive constants $C_1, C_2, C_3, C_4, \rho_1$ and $\rho_2$, such that, for any vector \textbf{v}, $C_1 \e1_{\|   \textbf{v}\| \leq \rho_1} \leq K_1(\textbf{v}) \leq C_2 \e1_{\|   \textbf{v}\| \leq \rho_1}$ and $C_3 \e1_{\|   \textbf{v}\| \leq \rho_1} \leq K_2(\textbf{v}) \leq C_4 \e1_{\|   \textbf{v}\| \leq \rho_2}$ . \\ \\
Condition (A8) is also imposed in the papers of Wei et al. (2013) and Xue (2009), where linear models with and without change-point are considered, respectively.\\
Concerning the selective probabilities functions $\pi_1(\textbf{x})$ and $\pi_2(\textbf{x})$, let us consider the following regularity hypotheses : \\  \\
\noindent \textbf{(A9)} $\pi_1(\textbf{x})$ and $\pi_2(\textbf{x})$ have bounded partial derivatives, with respect to $\textbf{x}$, up to order $\max (2, d-1)$ almost everywhere. \\ \\
Conditions (A7)-(A9) are usual assumptions for the convergence rate of the kernel estimation method, see for example the paper of Wei et al. (2013).
%----------------------------------------------------------------------
%--------------------------- Wighted Method ---------------------------
%----------------------------------------------------------------------
\subsubsection{Weighted empirical likelihood}
Now, we give the weighted empirical likelihood method. More specifically, in the vector $\textbf{g}$, we consider weight function of probabilities $\pi_1$, $\pi_2$. As discussed previously, complete-case analysis can yield biased estimates because the sample of observations that have no missing data might not be representative of the full sample. We could build a model to predict the nonresponse in that variable using all the other variables. The inverse of predicted probabilities of response from this model could then be used as survey weights to make the complete-case sample representative (along
the dimensions measured by the other predictors) of the full sample. This method becomes more complicated when there is more than one variable with missing data. \\

In order to obtain the weighted empirical likelihood statistic, we use the inverse probability weighted approach for missing data analysis, which was used by Horvitz and Thompson (1952) for missing data analysis. \\
We then define the two following $d-$random vectors
\begin{equation}
\label{eqW}
\left \{
\begin{array}{ccl}
\eg_{i,W}(\eb) \equiv \frac{\delta_i}{\pi_1(\XX_i)} \ef_i(\eb) (Y_i-f_i(\eb)),  &  &  i \in I, \\
\eg_{j,W}(\eb) \equiv \frac{\delta_j}{\pi_2(\XX_j)} \ef_j(\eb) (\tilde Y_j-f_j(\eb)),  &  &  j \in J.
\end{array}
\right.
\end{equation}
We recall that, if $Y_j$ non missing, $\tilde Y_j \equiv Y_j-f(\XX_j, \eb-\del^0)+f(\XX_j, \eb)$  for $j \in J$ and $\tilde Y_j$ equal to any finite value if $Y_j$ missing.\\ \\
Let us also consider in this case, the following matrices
\begin{equation*}
\label{VW}
\eV_{1n,W}(\eb) \equiv \frac{1}{k} \sum  _{i \in I}\egg_{i,W}(\eb), \qquad \eV_{2n,W}(\eb)    \equiv \frac{1}{n-k} \sum   _{j \in J} \egg_{j,W}(\eb).
\end{equation*}
\noindent Like as in the complete-case data, using similar argument of Remark \ref{rem1}, the test statistic for the weighted method is
\begin{equation}
\label{zWfin}
{\cal Z}_{nk,W}(\del_0,\el_W,\eb) \equiv 2 \Big[  \sum \limits _{i \in I} \log(1+\frac{k}{n} \el_W^t\et_{i,W}(\eb))+\sum \limits _{j \in J} \log(1-\frac{n}{n-k} \el_W^t \et_{j,W}(\eb) \Big],
\end{equation}
where
\begin{eqnarray*}
\et_{i,W}(\eb) \equiv \em_{n,W}^{\frac{1}{2}}(\eb)\eV_{1n,W}^{-1}(\eb)\eg_{i,W}(\eb), \qquad \et_{j,W}(\eb) \equiv \em_{n,W}^{\frac{1}{2}}(\eb)\eV_{2n,W}^{-1}(\eb)\eg_{j,W}(\eb),
\end{eqnarray*} 
\begin{eqnarray*}
\em_{n,W}(\eb)\equiv \frac{k(n-k)}{n^2} \eV_{1n,W}(\eb) \eh_{n,W}(\eb) \eV_{2n,W}(\eb)
\end{eqnarray*}
and
\begin{eqnarray*}
&&\eh_{n,W} (\eb) \equiv \Big [ \frac{ n-k}{nk} \sigma^2_1 \eV_{2n,W}(\eb) + \frac{k}{n(n-k)} \sigma^2_2 \eV_{1n,W}(\eb) \Big ]^{-1}. 
\end{eqnarray*}
The variances $ \sigma^2_1$ and $ \sigma^2_2$ can be estimated respectively by $\tilde \sigma^2_1$ and $\tilde \sigma^2_2$, given by (\ref{sigtild}). \\ \\
The score function of test statistic (\ref{zWfin}) is
\begin{equation}
\label{scoreZw}
\left\{
\begin{array}{ccl} 
{\ephi}_{1n,W}(\del_0,\el_W,\eb) = \displaystyle \sum \limits _{i \in I}
\frac{\et_{i,W}(\eb)}{\frac{k}{n}+\el^t_W \et_{i,W}(\eb)}-\sum \limits _{j\in J}
\frac{\et_{j,W}(\eb)}{\frac{n-k}{n}-\el^t_W\et_{j,W}(\eb)}, \\ \\
{\ephi}_{2n,W}(\del_0,\el_W,\eb)   = \displaystyle \sum \limits _{i \in I}
\frac{\ett_{i,W}(\eb)}{\frac{k}{n}+\el^t_W \et_{i,W}(\eb)}-\sum \limits _{j\in J}
\frac{\ett_{j,W}(\eb)}{\frac{n-k}{n}-\el^t_W\et_{j,W}(\eb)},
\end{array}
\right.
\end{equation} \ \\
where $\ett_{i,W}(\eb)$ and $\ett_{j,W}(\eb)$ are the derivative with respect to $\eb$ of $\et_{i,W}(\eb)$ and $\et_{j,W}(\eb)$ respectively.\\
Then, solving the system ${\ephi}_{1n,W}(\del_0,\el_W,\eb)=\textbf{0}_d$ and ${\ephi}_{2n,W}(\del_0,\el_W,\eb)=\textbf{0}_d$ given by (\ref{scoreZw}), we obtain $\hat \el_{n,W}(k)$ and $\hat \eb_{n,W}(k) $ the maximizers of the statistic (\ref{zWfin}). 
%----------------------------------------------------------------------
%--------------------------- Imputed Method ---------------------------
%----------------------------------------------------------------------
\subsubsection{Empirical likelihood with imputed values}
\hh In this part, we first reconstruct the missing response variable and then we propose the corresponding empirical likelihood. For the profile empirical likelihood with complete-case data and the weighted empirical likelihood, the information contained in the data is not explored fully. Since incomplete-case data are discarded in constructing the empirical likelihood ratio, the coverage accuracies of confidence regions are reduced when there are plenty of missing values. To resolve the issue, we use nonlinear regression imputation to impute $Y_i$( or $Y_j$) if $Y_i$( or $Y_j$) is missing. We introduce the forecast of $Y_l$, for $ l=1, \ldots, n$, constructed using the least square estimators for the parameters $\eb$ and $\eb_1$ and a nonparametric estimators for probabilities $\pi_1(\XX_i)$ and $\pi_2(\XX_j)$,
\begin{equation}
\label{eqq22}
\begin{array}{ccl}
 {Y}_{i,R} \equiv  \frac{ \delta_i}{\hat \pi_1(\XX_i)} Y_i+  \Big (1-\frac{ \delta_i}{\hat \pi_1(\XX_i)} \Big) f(\XX_i, \hat \eb) , &  & i\in I,  \\ \\
{Y}_{j,R} \equiv  \frac{ \delta_j}{\hat \pi_2(\XX_j)} Y_j+  \Big(1-\frac{ \delta_j}{\hat \pi_2(\XX_j)}\Big) f(\XX_j, \hat \eb_1 ),  &  & j \in J,
\end{array}
\end{equation}
where $\hat \eb$ and $\hat \eb_1$ given by relation (\ref{eqbb}).
\ \\ \ \\
We will consider the empirical likelihood on all these reconstructed (imputed) values. Then, for $ j \in J$, let us consider ${Y}_{j,R}^* \equiv {Y}_{j,R}-f(\XX_j, \eb-\del^0)+f(\XX_j, \eb)$. In this case, the auxiliary random vectors are defined by
\begin{equation*}
\label{eqI}
\left \{
\begin{array}{ccl}
\eg_{i,R}(\eb) \equiv  \ef_i(\eb) ({Y}_{i,R}-f_i(\eb)), &  &  i \in I, \\
\eg_{j,R}(\eb) \equiv \ef_j(\eb) ({Y}_{j,R}^*-f_j(\eb)), &  &  j \in J.
\end{array}
\right.
\end{equation*}
Let also the following matrices
\begin{equation*}
\label{VI}
\eV_{1n,R}(\eb) \equiv \frac{1}{k} \sum  _{i \in I} \egg_{i,R}(\eb), \qquad \eV_{2n,R}(\eb)    \equiv \frac{1}{n-k} \sum   _{j \in J}  \egg_{j,R}(\eb).
\end{equation*}
Like as in the two above methods, the test statistic for the imputed method is
\begin{equation}
\label{zIfin}
{\cal Z}_{nk,R}(\del_0,\el_{R},\eb) \equiv 2 \Big[  \sum \limits _{i \in I} \log(1+\frac{n}{k} \el_{R}^t\et_{i,R}(\eb))+\sum \limits _{j \in J} \log(1-\frac{n}{n-k} \el_{R}^t \et_{j,R}(\eb) \Big],
\end{equation}
where
\begin{equation*}
\et_{i,R}(\eb) \equiv \em_{n,R}^{\frac{1}{2}}(\eb) \eV_{1n,R}^{-1}(\eb)\eg_{i,R}(\eb), \qquad \et_{j,R}(\eb) \equiv \em_{n,R}^{\frac{1}{2}}(\eb) \eV_{2n,R}^{-1}(\eb)\eg_{j,R}(\eb),
\end{equation*}
\begin{eqnarray*}
\em_{n,R}(\eb)  \equiv \frac{k(n-k)}{n^2} \eV_{1n,R}(\eb) \eh_{n,R}(\eb) \eV_{2n,R}(\eb)
\end{eqnarray*}
and
\begin{eqnarray*}
\eh_{n,R} (\eb) &\equiv & \Big [ \frac{ n-k}{nk} \sigma^2_1 \eV_{2n,R}(\eb) \eV_{1n,R}^{-1}(\eb)  \sum \limits _{i \in I} \frac{1}{\pi_1(\textbf{X} _i)} \big [\eff_{i}(\eb)({Y}_{i,R}-f_i(\eb))-\ef_i(\eb)\ef_i^t(\eb) \big]
\\ &&
+ \frac{k}{n(n-k)} \sigma^2_2  \sum \limits _{j \in J} \frac{1}{\pi_2(\textbf{X} _j)} \big[\eff_j(\eb) ({Y}_{j,R}^*-f_j(\eb))-\ef_j(\eb)\ef_j^t(\eb) \big]\eV_{2n,R}^{-1}(\eb) \eV_{1n,R}(\eb) \Big ]^{-1}. 
\end{eqnarray*}
The score functions of test statistic of (\ref{zIfin}) are
\begin{equation}
\label{scoreZi}
\left\{
\begin{array}{ccl}
{\ephi}_{1n,R}(\del_0,\el_{R},\eb)  = \displaystyle \sum \limits _{i \in I}
\frac{\et_{i,R}(\eb)}{\frac{k}{n}+\el^t_{R} \et_{i,R}(\eb)}-\sum \limits _{j\in J}
\frac{\et_{j,R}(\eb)}{\frac{n-k}{n}-\el^t_{R}\et_{j,R}(\eb)}, \\ \\
{\ephi}_{2n,R}(\del_0,\el_{R},\eb)   =  \displaystyle \sum  \limits _{i \in I}
\frac{\ett_{i,R}(\eb)}{\frac{k}{n}+\el^t_{R} \et_{i,R}(\eb)}-\sum \limits _{j\in J}
\frac{\ett_{j,R}(\eb)}{\frac{n-k}{n}-\el^t_{R}\et_{j,R}(\eb)},
\end{array}
\right. 
\end{equation}
where $\ett_{i,R}(\eb)$ and $\ett_{j,R}(\eb)$ are the derivative with respect to $\eb$ of $\et_{i,R}(\eb)$ and $\et_{j,R}(\eb)$ respectively.\\
Then, solving the system ${\ephi}_{1n,R}(\del_0,\el_{R},\eb)=\textbf{0}_d$ and ${\ephi}_{2n,R}(\del_0,\el_{R},\eb)=\textbf{0}_d$ given by (\ref{scoreZi}), we obtain $\hat \el_{n,R}(k)$ and $\hat \eb_{n,R}(k)$ the maximizers of the statistic (\ref{zIfin}).
%%%%%%%%%%%%%%%%%%%%%%%%%%%%%%%%%%%%%%%%%%%%%
%%%%%%%%%%%%%%%%%%%%%%%%%%%%%%%%%%%%%%%%%%%%%
\subsection{Asymptotic behaviours of ${\cal Z}_{nk,C},{\cal Z}_{nk,W}$ and ${\cal Z}_{nk,R}$ }
\hspace*{0.3cm} In this subsection, we study the asymptotic distributions of the empirical likelihood ratios  ${\cal Z}_{nk,C},{\cal Z}_{nk,W}$ and ${\cal Z}_{nk,R}$, given by (\ref{zcfin}), (\ref{zWfin}) and (\ref{zIfin}) respectively. The main result is given by Theorem \ref{theo2}, where we show that under the null hypothesis $H_0$, all three statistics have, asymptotically, chi squared distributions. \\

We require the equivalent to the assumption (A5) given in the no missing response data case : \\ \\
\noindent  \textbf{(A10)} The matrices $\eV_{1n,C}$, $\eV_{2n,C}$, $\eV_{1n,W}$, $\eV_{2n,W}$, $\eV_{1n,R}$ and $\eV_{2n,R}$, are non singular for any $\XX_i\in \Upsilon$ and for all $\eb$ in a neighborhood of $\eb^0$. Moreover, the determinants of these matrices are bounded for sufficiently large n. \ \\ 
%Assumption (A10) assures that the matrices $\eV_{1n,C}$, $\eV_{2n,C}$, $\eV_{1n,W}$, $\eV_{2n,W}$, $\eV_{1n,Im}$ and $\eV_{2n,Im}$, are uniformly nonsingular and bounded for $n$ larger than some integer. \\

The following Proposition is similar to Proposition \ref{prop1}, for each of the three empirical likelihood ratio statistics with missing response variable, defined in subsection 4.1. So, we show that  $\|\hat \el_{n,C}(k)\|_2 \rightarrow 0$, $\|\hat \el_{n,W}(k)\|_2 \rightarrow 0$ and $\|\hat \el_{n,R}(k)\|_2 \rightarrow 0$, as $n\rightarrow \infty$ and that $\hat \eb_{n,C}(k)$, $\hat \eb_{n,W}(k)$ and $\hat \eb_{n,R}(k)$ are a consistent estimators of $\eb^0$, under hypothesis $H_0$. The proof is in appendix, Section 6.
%%%%%%%%%%%%%%%%%%%%%%%%%%%%%%%%%%%%%%%%%%%%%
%%%%%%%%%%%%%%%%%%%%%%%%%%%%%%%%%%%%%%%%%%%%%
\begin{prp}
\label{prop2}
Under null hypothesis $H_0$, if the assumptions (A1)-(A4), (A6)-(A10) are satisfied, then, for the estimators $(\hat \el_{n,C}(k),\hat \eb_{n,C}(k))$, $(\hat \el_{n,W}(k),\hat \eb_{n,W}(k))$ and $(\hat \el_{n,R}(k),\hat \eb_{n,R}(k))$ given by solving the systems (\ref{scoreZc}),(\ref{scoreZw}) and  (\ref{scoreZi}) respectively, we have 
\begin{eqnarray*}
\hat \el_{n,C}(k)=O_{\eP}(n^{-1/2}), \qquad  \hat \el_{n,W}(k)=O_{\eP}(n^{-1/2}), \qquad \hat \el_{n,R}(k) =O_{\eP}(n^{-1/2}), \\
\hat \eb_{n,C}(k)-\ebo=o_{\eP}(1), \qquad \hat \eb_{n,W}(k) -\ebo=o_{\eP}(1), \qquad \hat \eb_{n,R}(k)-\ebo=o_{\eP}(1).
\end{eqnarray*}
\end{prp}
\hh The following theorem gives the asymptotic distribution of the three empirical likelihood statistics given by (\ref{zcfin}), (\ref{zWfin}) and (\ref{zIfin}). The proof is in appendix, Section 6. 
\begin{thm}
\label{theo2}
Suppose that assumptions (A1)-(A4), (A6)-(A10) hold. Under null hypothesis $H_0$, the statistics ${\cal Z}_{nk,C}(\del_0,\el_C,\eb)$, ${\cal Z}_{nk,W}(\del_0,\el_W,\eb)$ and ${\cal Z}_{nk,R}(\del_0,\el_{R},\eb)$ all have an asymptotic $\chi^2_d$ distribution, where $d$ is the dimension of $\eb$. 
\end{thm}
Thus, for a fixed size $\alpha \in (0,1)$, an asymptotic $(1-\alpha)$ confidence region for $\del$, based on the empirical likelihood statistic for the three proposed methods are respectively,
\begin{eqnarray*}
&& CR_{\alpha,C}=\{ \del \in \R^d : {\cal Z} _ {nk,C} (\del_0,\el_C,\eb) < c_{1-\alpha;d} \},
\\ && 
CR_{\alpha,W}= \{ \del \in \R^d : {\cal {Z}}_ {nk,W}(\del_0,\el_W,\eb)  <c_{1-\alpha;d} \},
\\ &&
CR_{\alpha,R}=\{ \del \in \R^d :{\cal Z}_ {nk,R} (\del_0,\el_{R},\eb)  <c_{1-\alpha;d} \},
\end{eqnarray*}
where $c_{1-\alpha;d}$ is the $(1-\alpha)$ quantile of the chi-squared distribution with $d$ degrees of freedom.
 %%%%%%%%%%%%%%%%%%%%%%%%%%%%%%%%%%%%%%%%%%%%%%%%%
 %%%%%%%%%%%%%%%%%%%%%%%%%%%%%%%%%%%%%%%%%
 
\section{Simulation study}
In this subsection, we carried out some simulation studies to evaluate the performance of the proposed empirical likelihood confidence regions, by using Monte Carlo method. Firstly, when the response variable $Y_i$ is observed for each
observation $i=1, \ldots,n$, secondly, when some values of $\textbf{Y}$ may be missing. The program codes are available from the author.\ \\ 
%
%when the nonlinear regression model with complete-data, secondly, when we have a nonlinear regression model with missing response data. We use the R language for all simulations. 

We consider the following nonlinear function
\begin{equation} 
\label{model1}
f(x,\eb)=a\frac{1-x^b}{b},
\end{equation}
with $\eb=(a,b) \in [-100, 100] \times [0.1, 20]$. The same function was considered in Ciuperca (2011) where the model was estimated by penalized least absolute method and in Ciuperca and Salloum (2013) to test the change in the regression parameters of the nonlinear model.
\ \\
There are three different cases for the error terms :\ \\ \ \\
\textbf{Case a} : $\varepsilon_i ={\cal N}(0,1)$ and $\varepsilon_j ={\cal N}(0,1)$,\ \\
\textbf{Case b} : $\varepsilon_i =  1/\sqrt{6}(\chi^2(3)-3)$ and $\varepsilon_j=2/\sqrt{6}   t(6)$, 
\ \\
\textbf{Case c} : $\varepsilon_i=2{\cal E}xp(2)-1$ and $\varepsilon_j={\cal N}(0,1)$,\ \\ \ \\
where  ${\cal N}(0,1)$,  ${\cal E}xp(2)$, $\chi^2(3)$ and $t(6)$ are standard normal distribution, exponential distribution with  mean 1/2, chi-square distribution with degree of freedom 3 and Student distribution with degree of freedom 6, respectively.
% The nominal coverage level is $1-\alpha=0.95$.
 \ \\ \ \\
The following two nonlinear models are considered : \\ \\
 \textit{\textbf{Model 1}.} We generate the data sets from the following model
 \begin{equation}
\label{model1} 
\left\{ 
\begin{array}{ccl}
Y_i=a_1^0\frac{1-X_i^{b_1^0}}{b_1^0}+\varepsilon_i, &  & i=1, \ldots, k, \\
Y_j=a_1^0\frac{1-X_i^{b_1^0}}{b_1^0}+\varepsilon_j, &  & j=k+1, \ldots, n,
\end{array}
\right.
\end{equation}
where $a_1^0=10$, $b_1^0=2$ and $X_i=i/1000$. In this example,the null hypothesis is $H_0:\gamma_0=(0,0)^t$.\\ \\
 \textit{\textbf{Model 2}.} We generate the data sets from the following model
 \begin{equation}
\label{model2} 
\left\{ 
\begin{array}{ccl}
Y_i=a_1^0\frac{1-X_i^{b_1^0}}{b_1^0}+\varepsilon_i, &  & i=1, \ldots, k, \\
Y_j=a_2^0\frac{1-X_i^{b_2^0}}{b_2^0}+\varepsilon_j, &  & j=k+1, \ldots, n,
\end{array}
\right.
\end{equation}
where $a_1^0=10$, $b_1^0=2$, $a_2^0=7$, $b_2^0=1.75$ and $X_i=i/1000$. In this example, the null hypothesis is $H_0:\gamma_0=(3,0.25)^t$.

\subsection{\textbf{Model with complete data}}
\hspace*{0.3cm} For nominal confidence level $1- \alpha =0.95$ and 1000 Monte Carlo replications for each model/case, in Table 5.1 for the two models, we present the coverage probabilities (CP) and lengths of the confidence regions (LCR) obtained by the empirical log-likelihood method on the no-missing case data. These results were obtained for fixed $n=1000$ and different positions of change-point, using Theorem \ref{theo1}. Then, the asymptotic confidence region is $CR_{\alpha}=\{\del \in \R ^d :{\cal Z}_{nk}(\del_0,\el,\eb)  < c_{1-\alpha;d} \}$, where $c_{1-\alpha;d}$ is the $(1-\alpha)$ quantile of the standard chi-squared distribution with $d$ degrees of freedom, with ${\cal Z}_{nk}(\del_0,\el,\eb)$ given by equation (\ref{zfin}). \\
In order to calculate the coverage probability (CP), we consider a model under hypothesis $H_0$ given by (\ref{H0}) and we count the number of times, on the Monte Carlo replications, when the statistic value does not exceeds the critical value $c_{1-\alpha;d}$. The lengths of the confidence regions (LCR), for each model, designate the difference between value that we are confident of with upper or lower endpoint obtained for the statistic ${\cal Z}_{nk}$. 
We can see that all the coverage probabilities (CP) are very close to 0.95, which indicate the performance of the proposed empirical likelihood confidence region. The results don't varry with the point location $k$ or with the error distribution.
\begin{table}[h]
\label{tab1}
\begin{center}
\caption{Coverage probabilities (CP) and interval lengths (LCR) for models 1 and 2, with error distributions a, b, c.} 
%\ \\
\begin{tabular}{c c c c c c c c c c c c c}
   &    & \multicolumn{5}{c}{Model 1} &  & \multicolumn{5}{c}{ Model 2}\\ 
    \cline{3-7}  \cline{9-13}
%  &    & \multicolumn{5}{c}{ Model 1} &  & \multicolumn{5}{c}{ Model 2}\\ 
%    \cline{3-12}   
    $k$  & & \footnotesize Case a &  & \footnotesize
 Case b  &  & \footnotesize
 Case c & & \footnotesize
 Case a &  & \footnotesize
 Case b  &  & \footnotesize
 Case c \\ \noalign{\smallskip}\hline\noalign{\smallskip}
    300& CP &0.942& &0.933 & &0.949  & & 0.949 &  & 0.967  &  & 0.951\\
        & LCR & 4.876 & &5.682 & &5.470 & & 5.887&  & 5.978  &  & 5.957\\
    %\hline
 &   & & & & & & & &  &  &  & \\
   500 &  CP &0.966 & &0.900  & &0.960 & &0.958 &  &0.941  &  & 0.942\\
    & LCR  & 5.151 & &5.233 & & 4.939 & & 5.965 &  & 5.912  &  & 5.975 \\
    %\hline
    &   & & & & & & & &  &  &  & \\
    700 & CP &0.957 & &0.951  & &0.966 & & 0.951&  & 0.948  &  & 0.956\\
        & LCR &5.757 & &5.531 & &5.794 & & 5.976&  & 5.983  &  & 5.948 \\
    \hline
\end{tabular}
\end{center}
\end{table}
%%%%%%%%%%%%%%%%%%%%%%%%%%%%%%%%%
\subsection{Model with missing response data } 
\hspace*{0.3cm} In this subsection, we suppose that the response variable $Y_i$ can be missing at random.  \\
Throughout this subsection, the Kernel functions are taken as the Epanechnikov Kernel,\\ $K_1 (X_i)=0.75(1-X_i^2) \e1_{|X_i| \leq 1}$, $K_2 (X_j)=0.75(1-X_j^2) \e1_{|X_j| \leq 1}$ and the bandwidths $h_{1n}=k ^{-1/7}$, $h_{2n}=(n- k) ^{-1/7}$, which satisfy the condition (A7). \\
 \\
We consider the three following studies of response probabilities under the MAR assumption :  \\ \\
\textbf{Study 1} : $\pi_1(X_i)=0.8+0.2|X_i-1|$ if $|X_i-1|\leq 1$ and $0.95$ otherwise, $\pi_2(X_j)=0.8+0.2|X_j-1|$ if $|X_j-1|\leq 1$ and $0.95$ otherwise.
\\ \\
\textbf{Study 2} : $\pi_1(X_i)=0.8$ for all $X_i$ and $\pi_2(X_j)=0.8$ for all $X_j$.
\\ \\
\textbf{Study 3} : $\pi_1(X_i)=0.8+0.2|X_i-1|$ if $|X_i-1|\leq 1$ and $0.95$ otherwise, $\pi_2(X_j)=0.8$ for all $X_j$. \\ \\
For the studies 1, 2 and 3, the Tables 5.2, 5.3, and 5.4, present (CP) and (LCR) using statistics ${\cal Z}_{nk,C},{\cal Z}_{nk,W}$ and ${\cal Z}_{nk,R}$ given by equations (\ref{zcfin}), (\ref{zWfin}) and (\ref{zIfin}), respectively. For each study, $1-\alpha=0.95$, $n=1000$ and the error distributions are the same as in subsection 5.1. We run $1000$ replications for each simulation.
\\ \\
For Model 1, we obtain that the coverage probabilities (CP) are larger than 0.91. The results were very slightly lower of complete data case. On the other hand, the CP by weighted and imputed methods are the same as for complete-case. Conversely, the LCR are slightly wider for weighted and imputed methods. \\
For Model 2, we obtain the same results by the three methods.
\begin{table}[h]
\label{study1}       
\begin{center}
\caption{Coverage probability (CP) and interval lengths (LCR), of Study 1, $n=1000$, $k=600$.}
%\ \\
\begin{tabular}{c c c c c c c c c c c c c}
   &    & \multicolumn{5}{c}{Model 1} &  & \multicolumn{5}{c}{ Model 2} \\ 
    \cline{3-7} \cline{9-13}
   Y   & & \footnotesize Case a &  & \footnotesize
 Case b  &  & \footnotesize
 Case c & & \footnotesize
 Case a &  & \footnotesize
 Case b  &  & \footnotesize
 Case c    \\ \noalign{\smallskip}\hline\noalign{\smallskip}
     Comple-Case& CP &0.920  && 0.918 & & 0.913 && 0.956 && 0.937 && 0.951 \\
        & LCR & 4.450 &  & 5.249 &  & 4.988 && 5.964 && 5.944 && 5.953\\ 
          %  \noalign{\smallskip} \hline
            &   & & & & & & & &  &  &  & \\
Weighted &  CP & 0.917 &  & 0.910  &  & 0.903 && 0.948 && 0.938 && 0.949\\
    & LCR  & 5.896 &  & 5.792 &  & 4.725 && 5.982 && 5.988 && 5.914 \\
          %  \noalign{\smallskip}\hline
          &   & & & & & & & &  &  &  & \\
 Imputed & CP & 0.921 &  & 0.926  & &  0.915  && 0.951 && 0.942 && 0.953 \\
        & LCR & 5.630 &  & 5.759 &  & 4.778 && 5.857 && 5.979 && 5.915 \\
       \noalign{\smallskip}\hline
\end{tabular}
\end{center}
\end{table}
%%%%%%%%%%%%%%%%%%%%%%%%%%ù
%%%%%%%%%%%%%%%%%%%%%%%%%
\begin{table}[h]
\label{study2}       
\begin{center}
\caption{Coverage probability (CP) and interval lengths (LCR), of Study 2, $n=1000$, $k=600$.}
%\ \\
\begin{tabular}{c c c c c c c c c c c c c}
   &    & \multicolumn{5}{c}{Model 1} &  & \multicolumn{5}{c}{ Model 2} \\ 
    \cline{3-7} \cline{9-13}
   Y   & & \footnotesize Case a&  & \footnotesize
 Case b  &  & \footnotesize
 Case c & & \footnotesize
 Case a &  & \footnotesize
 Case b  &  & \footnotesize
 Case c    \\ \noalign{\smallskip}\hline\noalign{\smallskip}
     Comple-Case& CP &0.943 & & 0.911& & 0.932   && 0.951 && 0.943&& 0.932 \\
        & LCR & 5.963 & & 5.837 & &  5.869    && 5.947 && 5.973 && 5.971\\ 
          %  \noalign{\smallskip} \hline
            &   & & & & & & & &  &  &  & \\
Weighted &  CP & 0.901& & 0.905&  &0.918    && 0.967 && 0.942 && 0.938\\
    & LCR  & 5.981 & & 4.717 & & 4.778    && 5.906 && 5.959 && 5.967 \\
          %  \noalign{\smallskip}\hline
          &   & & & & & & & &  &  &  & \\
 Imputed & CP & 0.934& & 0.917& &0.919    && 0.954 && 0.946 && 0.952 \\
        & LCR & 5.566 & & 5.809 & & 4.594   && 5.964 && 5.926 && 5.964 \\
       \noalign{\smallskip}\hline
\end{tabular}
\end{center}
\end{table}
%%%%%%%%%%%%%%%%%%%%%%%%%%%
%%%%%%%%%%%%%%%%%%%
\begin{table}[h]
\label{study3}       
\begin{center}
\caption{Coverage probability (CP) and interval lengths (LCR), of Study 3, $n=1000$, $k=600$.}
%\ \\
\begin{tabular}{c c c c c c c c c c c c c}
   &    & \multicolumn{5}{c}{Model 1} &  & \multicolumn{5}{c}{ Model 2} \\ 
    \cline{3-7} \cline{9-13}
   Y   & & \footnotesize Case a &  & \footnotesize
 Case b  &  & \footnotesize
 Case c & & \footnotesize
 Case a &  & \footnotesize
 Case b  &  & \footnotesize
 Case c    \\ \noalign{\smallskip}\hline\noalign{\smallskip}
     Comple-Case& CP &0.925 & & 0.943 & & 0.914   && 0.950 && 0.957&& 0.948 \\
        & LCR & 5.676 & & 5.342 & & 5.863    && 5.982 && 5.976 && 5.982\\ 
          %  \noalign{\smallskip} \hline
            &   & & & & & & & &  &  &  & \\
Weighted &  CP & 0.902 & & 0.924 & & 0.913    && 0.951 && 0.932 && 0.924\\
    & LCR  & 4.854 & & 5.600 & & 5.376    && 5.965 && 5.839 && 5.920 \\
          %  \noalign{\smallskip}\hline
          &   & & & & & & & &  &  &  & \\
 Imputed & CP &  0.922 & & 0.928 & & 0.939    && 0.960 && 0.934 && 0.951 \\
        & LCR & 5.912 & & 5.559 & & 5.532   && 5.957 && 5.966&& 5.975 \\
       \noalign{\smallskip}\hline
\end{tabular}
\end{center}
\end{table}
%%%%%%%%%%%%%%%%%%%%
\ \\ \ \\ \ \\
\newpage
%\newpage
%%%%%%%%%%%%%%%%%%%%%%%%%%%%%%%%%%%%%%%
\section{Appendix}
We first give lemmas and their proofs, which are useful to prove propositions, theorems and of other lemmas. Then, we present the proofs of results stated  in Sections 3 and 4.
\subsection{Lemmas}
\begin{lem}
\label{lem00}
\noindent Let $\XX=(X_1,\ldots,X_p)$ a random vector (column), with the random variables $X_1,\ldots,X_p$ not necessarily independent and $\textbf{M}=(m_{uv})_{1 \leq u,v \leq p}$ such that $\textbf{M}= \XX \XX^t$. If for v=1, ..., p, we have
  \begin{equation}
\label{eqeta}
for\,\,all \,\,  \eta_v>0, there\,\,exists \,\, \delta_v >0\,\,\,\,\, such\,\, that \,\,\,\, \eP [|X_v|\geq \delta_v] \leq \eta_v, 
\end{equation}
then \ \\ 
 $(i)\,\, \eP \big[ \| \XX \|_1 \geq p \max  _{1 \leq v \leq p}\delta_v \big]\leq \max  _{1 \leq v \leq p}  \eta_v$, \\
 $(ii)\,\, \eP \big[ \| \XX \|_2 \geq \sqrt{p} \max  _{1 \leq v \leq p}\delta_v\big]\leq \max   _{1 \leq v \leq p}\eta_v$,  \\
 $(iii) \,\, \eP \big[ \| \textbf{M} \|_1 \geq p \max   _{1 \leq u,v \leq p} \{\delta_u^2,\delta_v^2\}\big]\leq \max   _{1 \leq u,v \leq p} \{\eta_u^2,\eta_v^2\}$, \ \\
 where $\| \textbf{M} \|_1= \max  _{1 \leq v \leq p} \{ \sum  \limits  _{u=1}^p |m_{uv}| \}$ is the subordinate norm to the vector norm $\|.\|_1$. 
 \end{lem}
\noindent {\bf Proof.} The proof of this lemma is given by Ciuperca and Salloum (2013). \hspace*{\fill}$\blacksquare$ 
%%%%%%%%%%%%%%%%%%%%%%%%%%%%%%%%%%%%%%%%
%%%%%%%%%%%%%%%%%%%%%%%%%%%%%%%%%%%%%%%%%%%%
%%%%%%%%%%%%%%%%%%%%%%%%%%%%%%%%%%%%%%%%%%%
\begin{lem}
\label{lem0}
Let the $\eta$-neighborhood of $\ebo$, ${\cal V}_{\eta}(\ebo)= \{ \eb \in \Gamma; \| \eb-\ebo\|_2 \leq \eta \}$, with $\eta \rightarrow 0$. Then, under assumptions (A1)-(A4), for all $\eb \in {\cal V}_{\eta}(\ebo)$, we have
\begin{equation}
\label{ev1}
\eV_{1n}(\eb)=\eV^0_{1n}+o_{\eP}(\eb-\ebo).
\end{equation}
We recall that, $\eV_{1n}(\eb)= \displaystyle \frac{1}{k} \sum \limits  _{i \in I} \egg_i(\eb)$.
\end{lem}
\noindent {\bf Proof.} By Taylor's expansion up to order 2 of $\egg_i(\eb)$ at $\eb=\ebo$, we obtain 
\begin{eqnarray}
\label{V1}
 \eV_{1n}(\eb)& =& \bigg [\frac{1}{k} \sum   _{i\in I} \eff_i (\ebo) \varepsilon_i - \frac{1}{k} \sum   _{i\in I} \ef_i (\ebo) \ef_i ^t(\ebo) \bigg ] 
\nonumber \\ &&
+ \frac{1}{k} \sum   _{i\in I} \em_i \varepsilon_i -\frac{1}{k} \sum   _{i\in I}\eff_i (\ebo) \ef_i  (\eb^{(1)}) (\eb-\ebo)^t
\nonumber \\ &&
- \frac{1}{k} \sum   _{i\in I} \em_i  \ef_i  (\eb^{(1)})(\eb-\ebo)^t
-\frac{1}{k} \sum   _{i\in I} \em_{1i}  \ef_i  (\ebo) (\eb-\ebo)^t
\nonumber \\ &&
 + \frac{1}{k} \sum   _{i\in I}  \em_{2i} \ef_i  (\ebo)(\eb-\ebo)^t
+\frac{1}{k} \sum   _{i\in I}\em_{2i} \em_{3i}  (\eb-\ebo) (\eb-\ebo)^t, 
\end{eqnarray}
where, $\em_i$, $\em_{1i}$ and $ \em_{2i}$ are a $ d \times d$ square matrices, defined by 
$\em_i=  \sum _{r=1}^d \frac{\partial \eff_i (\eb_{i,r}^{(2)}) }{ \partial \beta_r } (\beta_r - \beta^0_r) $, $\em_{1i}= \big (\frac{\partial ^2 f_i (\eb_{i,rl}^{(3)}) }{ \partial \beta_r \partial \beta_l }  \big )_{1 \leq r,l \leq d}$ and 
$\em_{2i}=\big ( \frac{\partial ^2 f_i (\eb_{i,rl}^{(4)}) }{ \partial \beta_r \partial \beta_l }  \big )_{1 \leq r,l \leq d}$. Here, for $1 \leq r,l \leq d$, $ \beta_r$ denotes the r-th component of $\eb $ and $\eb^{(a)}_{i,rl}= \ebo+ u^{(a)}_{i,rl} (\eb-\ebo)$, with $u^{(a)}_{i,rl} \in [0,1]$ and $ a \in \{2,3,4\}$. \\ \\
For the first term of the right-hand side of (\ref{V1}), we have 
\begin{equation}
\label{V1,0}
\frac{1}{k} \sum   _{i\in I} \eff_i (\ebo)  \varepsilon_i - \frac{1}{k} \sum   _{i\in I} \ef_i (\ebo) \ef_i ^t(\ebo)= \eV_{1n}^0.
\end{equation}
By Bienaymé-Tchebychev's inequality and assumption (A1), we obtain that for all $C_1 > 0$ and $ i \in I$ 
 \begin{equation}
\label{epselon}
\eP \big[| \varepsilon_i | > C_1 \big] \leq  \frac{\sigma_1^2}{C_1}.
\end{equation}
For $ 1  \leq r,l \leq d $ and for any fixed $i$ such that $ i \in I$, denote by $M_{i,sl}$ the following random variable designating the term $(s,l)$ of the matrix $\em_i$ such that
$$M_{i,sl}=  \sum  _{r=1}^d \frac{\partial f_i(\eb_{i,r}^{(2)})}{\partial \beta_r \partial \beta_s \partial \beta_l } (\beta_r-\beta^0_r).$$
Using assumption (A3), we have with a probability one that $ |M_{i,sl}| \leq \| \eb - \ebo \|_2$. Applying Lemma \ref{lem00}(iii), for all $C_2 > 0$ and $ i \in I$, we obtain that 
\begin{equation}
\label{Mi}
\| \textbf{M}_i \|_1 \leq C_2 \| \eb - \ebo \|_2.
\end{equation}
For the second term of the right-hand side of (\ref{V1}), we have 
$$\|\frac{1}{k} \sum   _{i\in I} \em_i \varepsilon_i \|_1   \leq \frac{1}{k} \sum   _{i\in I} \| 
\em_i \varepsilon_i \|_1 \leq  \frac{1}{k} \sum   _{i\in I} \| 
\em_i \|_1 | \varepsilon_i | .$$
 Using relations (\ref{epselon}) and (\ref{Mi}), we obtain that
 \begin{equation}
\label{Rel1} 
\frac{1}{k} \sum   _{i\in I} \em_i \varepsilon_i = o_{\eP} ( \eb -\ebo).
\end{equation}
By Markov's inequality, taking also in account assumption (A4), then for any fixed $i$ and for $ 1 \leq r \leq l$, we have that $\eP[ |\frac{\partial f_i  (\eb_{i,r}^{(2)})}{ \partial \beta_r} | > \epsilon_1 ] \leq \eE(|\partial f_i  (\eb_{i,r}^{(2)})/\partial \beta_r | ) / \epsilon_1  $. For $\epsilon_1=\eE(|\partial f_i  (\eb_{i,r}^{(2)})/\partial \beta_r | ) / \epsilon  $, we obtain that $\eP[ | \partial f_i  (\eb_{i,r}^{(2)})/\partial \beta_r | > \eE(|\partial f_i  (\eb_{i,r}^{(2)})/\partial \beta_r |) / \epsilon   ] \leq \epsilon $. This last relation together with Lemma \ref{lem00}(i) imply
 \begin{equation}
\label{Rel2}
\eP[ \|\ef_{i} (\eb_{i,r}^{(2)}) \|_1 \geq \frac{d}{\epsilon} \max (\eE(|\partial f_i  (\eb_{i,r}^{(2)})/\partial \beta_r | )  ] \leq \epsilon.
\end{equation}
%In the same way, by supposition (A3) and lemma 1(i) of \cite {Ciuperca:Salloum}, for all $\epsilon_2 >0$, we have 
%\begin{equation}
%\label{Rel3}
%\eP[ \| \eff (\ebo)   \|_1 \geq \epsilon_2 ] \leq \epsilon.
%\end{equation}
For the third term of the right-hand side of (\ref{V1}), by assumption (A3) and using relation (\ref{Rel2}), we obtain that $k^{-1} \sum   _{i\in I} \eff_i ( \ebo) \ef_i  (\eb^{(1)}) (\eb-\ebo)^t\|_1 \leq (k^{-1} \sum   _{i\in I} \| \eff_i ( \ebo)  \ef_i  (\eb^{(1)})(\eb - \ebo)^t \|_1 \leq k^{-1} \sum   _{i\in I} \|\eff_i ( \ebo) \|_1$ $\\ \cdot  \| \ef_i  (\eb^{(1)}) \|_1 \|\eb - \ebo \|_1 $, which gives 
\begin{equation}
\label{Rel4}
\frac{1}{k} \sum   _{i\in I} \eff_i ( \ebo) \ef_i  (\eb^{(1)})(\eb-\ebo)^t= o_{\eP} ( \eb -\ebo).
\end{equation}
Using assumption (A3), we obtain that for $  1\leq s,r \leq d$, for all $\epsilon>0$ there exists $\epsilon_3>0$ such that $
\eP \big[|\frac{\partial^2 f_i( \eb^{(3)}_{i,sr})}{\partial \beta_s \partial \beta_r} | \geq  \epsilon_2\big] \leq \epsilon$.
By Lemma \ref{lem00} (iii), we have that for all $\epsilon>0$, 
\begin{equation}
\label{eqbt14}
\eP \big[\| \em_{1i}  \|_1 \geq  \epsilon_3 \big] \leq \epsilon.
\end{equation}
Using assumption (A3) and by a similar arguments as $\em_{1i}$, we can demonstrate that, for $  1\leq s,r \leq d$, for all $\epsilon>0$ there exists $\epsilon_4>0$ such that 
\begin{equation}
\label{eqbt15}
\eP \big[\| \em_{2i}  \|_1 \geq  \epsilon_4\big] \leq \epsilon.
\end{equation}
For the fourth term of the right-hand side of (\ref{V1}), using relations (\ref{Mi}),  (\ref{Rel2}) and the fact that $\| \eb - \ebo \| _2 \leq \eta $, we obtain that 
\begin{equation}
\label{Rel5}
\frac{1}{k} \sum   _{i\in I} \em_i \ef_i  (\eb^{(1)})(\eb-\ebo)^t= o_{\eP} ( \eb -\ebo).
\end{equation}
Using assumptions (A3) and relations (\ref{Rel2}), (\ref{eqbt14}), for the fifth term of the right-hand side of (\ref{V1}), we have
\begin{equation}
\label{Rel6}
\frac{1}{k} \sum   _{i\in I}  \em_{1i} \ef_i  (\ebo) (\eb-\ebo)^t= o_{\eP} ( \eb -\ebo).
\end{equation}
In the same way, using assumption (A3) and relations (\ref{Rel2}), (\ref{eqbt14}), (\ref{eqbt15}), we obtain that
\begin{equation}
\label{Rel7}
\frac{1}{k} \sum   _{i\in I} (\eb-\ebo) \em_{2i} \ef_i ^t (\ebo)= o_{\eP} ( \eb -\ebo)
\end{equation}
and
\begin{equation}
\label{Rel8}
\frac{1}{k} \sum   _{i\in I} (\eb-\ebo) \em_{2i} \em_{3i} (\eb-\ebo)^t = o_{\eP} ( \eb -\ebo).
\end{equation}
Combining relations (\ref{V1,0}), (\ref{Rel1}), (\ref{Rel4}) and (\ref{Rel5})-(\ref{Rel8}), we obtain that $\eV_{1n}(\eb)=\eV_{1n}^0+o_{\eP} ( \eb -\ebo)$. \hspace*{\fill}$\blacksquare$  \ \\  \ \\
Similarly to Lemma \ref{lem0}, we can demonstrate easily that $\eV_{2n}(\eb)=\eV_{2n}^0+o_{\eP} ( \eb -\ebo)$. Consequently, we have $\em_n(\eb)=\em_n(\ebo)+o_{\eP}(\eb-\ebo)$, with $\em_n(\eb)$ given in Remark \ref{rem1}. To simplify notation we write $\em_n(\ebo)=\em_n^0$.\ \\  \ \\
%%µµµµµµµµµµµµµµµµµµµµµµµµµµµµµµµµµµµµµµµµµµµµµµµµµµµµµµµµµµµµµµµ
%µµµµµµµµµµµµµµµµµµµµµµµµµµµµµµµµµµµµµµµµµµµµµµµµµµµµµµµµµµµµµµµ
%µµµµµµµµµµµµµµµµµµµµµµµµµµµµµµµµµµµµµµµµµµµµµµµµµµµµµµµµµµµµµµµ
%µµµµµµµµµµµµµµµµµµµµµµµµµµµµµµµµµµµµµµµµµµµµµµµµµµµµµµµµµµµµµµµ
%µµµµµµµµµµµµµµµµµµµµµµµµµµµµµµµµµµµµµµµµµµµµµµµµµµµµµµµµµµµµµµµ
The next lemma gives the behavior of $\sum  \limits _{i\in I}\et_i(\eb)$ in the neighborhood of $\ebo$.
\begin{lem} 
\label{lem1}
Let the $\eta$-neighborhood of $\ebo$, ${\cal V}_{\eta}(\ebo)= \{ \eb \in \Gamma; \| \eb-\ebo\|_2 \leq \eta \}$, with $\eta \rightarrow 0$. Then, under assumptions (A1)-(A6), for all $\eb \in {\cal V}_{\eta}(\ebo)$, we have
\begin{equation}
\label{zi0}
\frac{1}{k} \sum   _{i\in I}\et_i(\eb)=\frac{1}{k} \sum   _{i\in I}\et_i(\ebo)-\frac{1}{2} (\em_n^0)^{ \frac{1}{2}} (\eb-\ebo)+o_{\eP}(\eb-\ebo).
\end{equation}
\end{lem}
\noindent {\bf Proof.} Under assumption (A5), using Lemma \ref{lem0}, then the Taylor's expansion up to the order 2 of $k^{-1} \sum  _{i\in I}\et_i(\eb)$ at $\eb=\eb^0$ is
\begin{eqnarray}
\label{dlTi}
\frac{1}{k} \sum   _{i\in I}\et_i(\eb)&=& \frac{1}{k} \sum   _{i\in I}(\em_n^0)^{ \frac{1}{2}} (\eV_{1n}^0)^{-1}\ef_i(\ebo) \varepsilon_i
 +\frac{1}{2k} \sum   _{i\in I}(\em_n^0)^{ \frac{1}{2}} (\eV_{1n}^0)^{-1} \eff_i (\ebo) (\eb-\ebo) \varepsilon_i 
 \nonumber \\ &&
 -\frac{1}{2k} \sum   _{i\in I}(\em_n^0)^{ \frac{1}{2}} (\eV_{1n}^0)^{-1}\ef_i(\ebo)\ef_i^t (\ebo) (\eb-\ebo) \nonumber \\ && -\frac{1}{6k} \sum   _{i\in I}(\em_n^0)^{ \frac{1}{2}} (\eV_{1n}^0)^{-1}\ef_i(\ebo)(\eb-\ebo)^t \em_{4i} (\eb-\ebo) 
 \nonumber \\ &&
- \frac{1}{4k} \sum   _{i\in I} (\em_n^0)^{ \frac{1}{2}} (\eV_{1n}^0)^{-1} \eff_i (\ebo) (\eb-\ebo)\ef_i^t (\ebo) (\eb-\ebo)
\nonumber \\ &&
 -\frac{1}{12k} \sum   _{i\in I} (\em_n^0)^{ \frac{1}{2}} (\eV_{1n}^0)^{-1} \eff_i (\ebo) (\eb-\ebo)(\eb-\ebo)^t \em_{4i} (\eb-\ebo)
 \nonumber \\ &&
 - \frac{1}{12k} \sum   _{i\in I} (\em_n^0)^{ \frac{1}{2}} (\eV_{1n}^0)^{-1}\textbf{D}_i (\eb-\ebo)^t \em_{4i} (\eb-\ebo)
 \nonumber \\ &&
 + \frac{1}{6k} \sum   _{i\in I} (\em_n^0)^{ \frac{1}{2}} (\eV_{1n}^0)^{-1} \textbf{D}_i \varepsilon_i,
\end{eqnarray}
where $\em_{4i}= \big (\frac{\partial ^2 f_i (\eb_{i,ul}^{(5)}) }{ \partial \beta_u \partial \beta_l }  \big )_{1 \leq u,l \leq d}$ is a $d \times d $ matrix and 
$$\textbf{D}_i= \Big( \sum  _{l=1}^d \sum  _{u=1}^d \frac{\partial^2 \ef_i (\eb_{i,ul}^{(6)})}{\partial \beta_u \partial \beta_l} 
\cdot  (\beta_u-\beta^0_u)(\beta_l-\beta^0_l)\Big) $$
 is  a vector of dimension $(d \times 1)$. For $1 \leq u,l \leq d$, $ \beta_u$ denotes the u-th component of $\eb $ and $\eb^{(a)}_{i,ul}= \ebo+ v^{(a)}_{i,ul} (\eb-\ebo)$, with $v^{(a)}_{i,ul} \in [0,1]$ and $ a \in \{5,6 \}$. \\ \\ 
For the second term of the right-hand side of (\ref{dlTi}), by the law of large numbers, the term $(2k)^{-1}$ \\ $ \sum   _{i\in I}(\em_n^0)^{ \frac{1}{2}} (\eV_{1n}^0)^{-1}\eff_i(\ebo) (\eb-\ebo) \varepsilon_i $ converges almost surely to the expected of $(\em_n^0)^{ \frac{1}{2}} (\eV_{1n}^0)^{-1}\eff_i(\ebo) $ \\ $ \cdot(\eb-\ebo) \varepsilon_i $, as $n \rightarrow \infty $. Since $\varepsilon_i$ is independent of $\XX_i$ and $\eE[\varepsilon_i]=0$, we have
\begin{equation}
\label{d1}
\frac{1}{2k} \sum   _{i\in I}(\em_n^0)^{ \frac{1}{2}} (\eV_{1n}^0)^{-1}\eff_i(\ebo) (\eb-\ebo) \varepsilon_i=o_{\eP}(\eb-\ebo).
\end{equation}
For the third term of the right-hand side of (\ref{dlTi}), by a simple computation, we can obtain
\begin{equation}
\label{d2}
\frac{1}{2k} \sum   _{i\in I}(\em_n^0)^{ \frac{1}{2}} (\eV_{1n}^0)^{-1}\ef_i(\ebo)\ef_i^t (\ebo) (\eb-\ebo)=-\frac{1}{2} (\em_n^0)^{ \frac{1}{2}} (\eb-\ebo).
\end{equation}
Using assumption (A3), (A4) and by an similar arguments to those for relations (\ref{Rel2}) and (\ref{eqbt14}), we have for the fourth, fifth and the sixth terms of the right-hand side of (\ref{dlTi}), respectively, that
\begin{equation}
\label{d3}
\frac{1}{6k} \sum   _{i\in I}(\em_n^0)^{ \frac{1}{2}} (\eV_{1n}^0)^{-1}\ef_i(\ebo)(\eb-\ebo)^t \em_{4i} (\eb-\ebo)=o_{\eP}(\eb-\ebo),
\end{equation}
%In the same way, using assumption (A3) and Lemma 3 of \cite{Ciuperca:Salloum}, for the fifth and the sixth terms of the right-hand side of (\ref{dlTi}),  we have respectively
\begin{equation}
\label{d4}
\frac{1}{4k} \sum   _{i\in I} (\em_n^0)^{ \frac{1}{2}} (\eV_{1n}^0)^{-1}\eff_i(\ebo) (\eb-\ebo)\ef_i^t (\ebo) (\eb-\ebo)=o_{\eP}(\eb-\ebo)
\end{equation}
and
\begin{equation}
\label{d5}
\frac{1}{12k} \sum   _{i\in I} (\em_n^0)^{ \frac{1}{2}} (\eV_{1n}^0)^{-1} \eff_i(\ebo) (\eb-\ebo)(\eb-\ebo)^t \em_{4i} (\eb-\ebo)=o_{\eP}(\eb-\ebo).
\end{equation}
For any fixed $i$, such that $1 \leq i \leq k$ and for $1\leq s \leq d$, let us denote by $D_{is}$ the following random variable designating the s-th component of the vector $\textbf{D}_i$, such that
$$\textbf{D}_{is}= \sum  _{l=1}^d \sum  _{u=1}^d \frac{\partial^3 f_i(\eb_{i,ul}^{(6)})}{\partial \beta_u \partial \beta_l \partial \beta_s} (\beta_u-\beta^0_u)(\beta_l-\beta^0_l).$$
Applying Lemma \ref{lem00} and by assumption (A3), we have for all $C_3 >0$
\begin{equation}
\label{eqD}
\| \textbf{D} \|_1 \leq C_3 \| \eb-\ebo\|_2^2.
\end{equation}
The above equation, together with (\ref{epselon}), implies that
\begin{equation}
\label{d6}
\frac{1}{6k} \sum   _{i\in I} (\em_n^0)^{ \frac{1}{2}} (\eV_{1n}^0)^{-1} \textbf{D}_i \varepsilon_i =o_{\eP}(\eb-\ebo).
\end{equation}
Finally, for the term $(12k)^{-1} \sum   _{i\in I} (\em_n^0)^{ \frac{1}{2}} (\eV_{1n}^0)^{-1}\textbf{D}_i (\eb-\ebo)^t \em_{4i} (\eb-\ebo)$, assumption (A3), together with relation (\ref{eqD}) yield
\begin{equation}
\label{d7}
\frac{1}{12k} \sum   _{i\in I} (\em_n^0)^{ \frac{1}{2}} (\eV_{1n}^0)^{-1}\textbf{D}_i (\eb-\ebo)^t \em_{4i} (\eb-\ebo)=o_{\eP}(\eb-\ebo).
\end{equation}
Combining relations (\ref{d1})-(\ref{d5}), (\ref{d6}) and (\ref{d7}), lemma yields.\hspace*{\fill}$\blacksquare$  \ \\

%*************************************************************
%******************************************************************
%****************************************************************************
%**************************************************
 
%µµµµµµµµµµµµµµµµµµµµµµµµµµµµµµµµµµµµµµµµµµµµµµµµµµµµµµµµµµµµµµµ
%\label{proof prop 1}
\subsection{Proposition and Theorem proofs}
\noindent {\bf Proof of Proposition \ref{prop1}.} By the definition of the empirical likelihood ratio, we have the constraints $\sum  \limits  _{i\in I}p_i \eg_i( \eb)=\sum  \limits  _{j\in J}q_j \eg_j( \eb)=\textbf{0}_d$, which give $\sum \limits  _{i\in I}p_i \et_i( \eb)=\sum  \limits  _{j\in J}q_j \et_j( \eb)=\textbf{0}_d$.\\
Using the value of $p_i$ and $q_j$ given by (\ref{eqq1}) and (\ref{eqq2}) respectively and by an elementary calculations, we obtain
\begin{equation}
\label{pizi}
\sum   _{i\in I}p_i \et_i( \eb) =\frac{1}{k} \sum  _{i \in I} \et_i(\eb) -\frac{n}{k^2} \sum  _{i \in I} \frac{\et_i(\eb) \et_i^t (\eb) }{1+\frac{n}{k} \el^t(\eb)\et_i(\eb)}\el(\eb) = \textbf{0}_d
\end{equation}
and
\begin{equation}
\label{qjzj}
\sum  _{j\in J}q_j \et_j( \eb) =\frac{1}{n -k} \sum  _{j\in J} \et_j(\eb) +\frac{n}{(n-k)^2} \sum  _{j\in J} \frac{\et_j(\eb) \et_j^t (\eb) }{1+\frac{n}{n-k} \el^t(\eb)\et_j(\eb)}\el(\eb) = \textbf{0}_d.
\end{equation}
For the term $k^{-1} \sum   _{i \in I} \et_i(\eb)$ of (\ref{pizi}), by Lemma \ref{lem1} we have that 
\begin{equation}
\frac{1}{k} \sum   _{i\in I}\et_i(\eb)=\frac{1}{k} \sum   _{i\in I}\et_i(\ebo)-\frac{1}{2} (\em_n^0)^{ \frac{1}{2}} (\eb-\ebo)+o_{\eP}(\eb-\ebo).
\end{equation}
For the term $(n /k^2) \sum \limits  _{i \in I} \frac{\et_i(\eb) \et_i^t (\eb) }{1+\frac{n}{k} \el^t(\eb)\et_i(\eb)}\el(\eb)$, using Proposition 1 of Ciuperca and Salloum (2013), we have that for all $\epsilon >0$, there exists $N_1, N_2 > 0$, such that   
\begin{equation*}
\label{zizit}
\eP \Big [\frac{1}{N_1}  \sum _{i \in I} \et_i(\eb) \et^t_i(\eb)\leq  \sum   _{i \in I}\frac{\et_i(\eb) \et^t_i(\eb)}{1+\frac{n}{k} \el^t(\eb)\et_i(\eb)} \leq \frac{1}{N_2}  \sum  _{i \in I} \et_i(\eb) \et^t_i(\eb) \Big ]<\epsilon.
\end{equation*}
This implies that, in order to study the second term of the right-hand side of (\ref{pizi}), we must study only $k^{-1} \sum  _{i \in I} \et_i(\eb) \et^t_i(\eb)$.\ \\
By Lemma \ref{lem1}, we have that $\eV_{1n}(\eb)=\eV_{1n}^0+ o_{\eP} (\eb - \ebo)$. Under assumption (A5), the Taylor's expansion up to the order 2 of $\et_i(\eb)$ at $\eb=\eb^0$ is
\begin{eqnarray}
\label{devzi}
\et_i(\eb)&=& (\em_n^0)^{ \frac{1}{2}} (\eV_{1n}^0)^{-1}\ef_i(\ebo) \varepsilon_i+\frac{1}{2}(\em_n^0)^{ \frac{1}{2}} (\eV_{1n}^0)^{-1}\em_{6i} (\eb-\ebo) \varepsilon_i \nonumber  \\ &&  -\frac{1}{2}(\em_n^0)^{ \frac{1}{2}} (\eV_{1n}^0)^{-1}\ef_i(\ebo)\ef_i^t (\ebo) (\eb-\ebo) \nonumber  \\ &&
 -\frac{1}{6}(\em_n^0)^{ \frac{1}{2}} (\eV_{1n}^0)^{-1}\ef_i(\ebo)(\eb-\ebo)^t \em_{7i} (\eb-\ebo)  \\ &&  - \frac{1}{4}(\em_n^0)^{ \frac{1}{2}} (\eV_{1n}^0)^{-1}\em_{6i} (\eb-\ebo)\ef_i^t (\ebo) (\eb-\ebo) \nonumber\\ &&
-\frac{1}{12}(\em_n^0)^{ \frac{1}{2}} (\eV_{1n}^0)^{-1} \em_{6i}(\eb-\ebo)(\eb-\ebo)^t \em_{7i} (\eb-\ebo), \nonumber
\end{eqnarray}
where, $\em_{6i}$ and $ \em_{7i}$ are a $ d \times d$ square matrices, defined by 
$\em_{6i}= \big (\frac{\partial ^2 f(\eb_{i,rl}^{(7)}) }{ \partial \beta_r \partial \beta_l }  \big )_{1 \leq r,l \leq d}$, 
$\em_{7i}=\big ( \frac{\partial ^2 f(\eb_{i,rl}^{(8)}) }{ \partial \beta_r \partial \beta_l }  \big )_{1 \leq r,l \leq d}$  and $ \eb^{(a)}_{i,rl}=\ebo+v_{i,rl}^{(a)}(\eb- \ebo)$, with $v_{i,rl}\in [0,1]$ and $ a \in \{7,8 \}$. \\ \\
For the second term of (\ref{pizi}), using the Taylor's expansion of $\et_i(\eb)$ in a neighborhood of $\ebo$ given by the relation (\ref{devzi}) and with a similar argument to the one used in Lemma \ref{lem1} for the first term of (\ref{pizi}), together with the assumptions (A3), (A4), we obtain 
\begin{eqnarray}
\label{zi0zi0}
\frac{1}{k} \sum   _{i\in I}\et_i(\eb)\et^t_i(\eb)=\frac{1}{k} \sum   _{i\in I}\et_i(\ebo)\et^t_i(\ebo)+\eephi^0_{1n}+o_{\eP}( \eb-\ebo),
\end{eqnarray}
where \\
$\eephi^0_{1n}=\displaystyle \frac{n}{k^2} \sum   _{i\in I} (\em_n^0)^{ \frac{1}{2}} (\eV_{1n}^0)^{-1}\ef_i(\ebo)\ef_i^t(\ebo)(\eb-\ebo)  \big [(\em_n^0)^{ \frac{1}{2}} (\eV_{1n}^0)^{-1}  \ef_i(\ebo) \ef_i^t(\ebo)(\eb-\ebo) \big]^t$. \\ \\
In the same way, for the observations $ j \in J$, we obtain
\begin{equation}
\label{zj0}
\frac{1}{n-k} \sum   _{j\in J}\et_j(\eb)=\frac{1}{n-k} \sum  _ {j\in J}\et_j(\ebo)-\frac{1}{2} (\em_n^0)^{ \frac{1}{2}} (\eb-\ebo)+o_{\eP}(\eb-\ebo)\\
\end{equation}
and 
\begin{eqnarray}
\label{zj0zj0}
\frac{1}{n-k} \sum   _{j\in J}\et_i(\eb)\et^t_j(\eb) = \frac{1}{n-k} \sum   _{j\in J}\et_j(\ebo)\et^t_j(\ebo)+\eephi^0_{2n} +o_{\eP}( \eb-\ebo), \nonumber
\end{eqnarray}
where \\ 
$\eephi^0_{2n}= \displaystyle \frac{n}{(n-k)^2} \sum  \limits _{j\in J} (\em_n^0)^{ \frac{1}{2}} (\eV_{2n}^0)^{-1} \ef_j(\ebo)\ef_j^t(\ebo)(\eb-\ebo) \big [(\em_n^0)^{ \frac{1}{2}} (\eV_{2n}^0)^{-1} \ef_j(\ebo)\ef_j^t(\ebo)(\eb-\ebo) \big]^t$. \\ \\
To facilitate writing, we consider the $ d \times d$ square matrices, defined by \ \\ \ \\
$$\es_n^0=\frac{n}{k^2} \sum \limits  _{i\in I}\et_i(\ebo)\et^t_i(\ebo)+\frac{n}{(n-k)^2} \sum  \limits _{j\in J}\et_j(\ebo)\et^t_j(\ebo),$$ 
$$\eephi_n^0=\eephi_{1n}^0+\eephi^0_{2n}$$
and we define the vector \ \\ \ \\
$$\epsi_n ^0=\frac{1}{k} \sum  \limits _{i\in I}\et_i(\ebo)-\frac{1}{n-k} \sum \limits  _{j\in J}\et_j(\ebo).$$
\noindent On the other hand, we have $\ephi_{1n}(\del_0,\hat \el_n (k),\hat \eb_n (k))=\textbf{0}_d$. Using relations (\ref{zi0}) and (\ref{zi0zi0})-(\ref{zj0zj0}), we obtain
%%%%%%%%%%%%%%%%%%%%%%%%%%%
\begin{eqnarray*}
&& \Big[\frac{1}{k}  \sum \limits  _{i\in I}\et_i(\ebo)-\frac{1}{2} (\em_n^0)^{ \frac{1}{2}} (\hat \eb_n( k)-\ebo)- \eephi^0_{1n}\hat \el_n ( k)- \frac{n}{k^2} \sum \limits  _{i\in I} \et_i(\ebo) \et^t_i(\ebo) \hat \el_n( k)\Big]
\nonumber \\  &&
- \Big[\frac{1}{k} \sum \limits _ {j\in J} \et_j(\ebo)-\frac{1}{2} (\em_n^0)^{ \frac{1}{2}} (\hat \eb_n( k)-\ebo)+\eephi^0_{2n}\hat \el_n( k)+\frac{n}{(n-k)^2} \sum \limits  _{_j\in J} \et_j(\ebo) \et^t_j(\ebo)\hat \el_n( k)\Big]
\nonumber \\  &&
= \textbf{0}_d.
\end{eqnarray*}
Using notations given above, then we obtain
\begin{equation}
\label{el1}
\hat{\el}_n( k)=(\es_n^0+\eephi_n^0)^{-1} \epsi_n^0 + o_{\eP}(\hat \eb _n( k)-\ebo).
\end{equation}
The limited development of the statistic ${\cal Z}_{nk}(\del_0,\hat \el_n( k),\hat\eb_n( k))$ specified by relation (\ref{zfin}), in the neighborhood of $(\el,\eb)=(\textbf{0}_d,\ebo)$ up to the order 2, can be written
\begin{eqnarray}
\label{DLz}
&& {\cal Z}_{nk}(\del_0,\hat \el_n( k),\hat \eb _n( k))=\Big [2\hat {\el}^t _n( k)\Big( \frac{n}{k}  \sum  _{i \in I}\et_i(\ebo)-  \frac{n}{n-k} \sum  _{j \in J} \et_j(\ebo)\Big) \Big]  \nonumber \\ &&
- \Big [ \hat {\el}^t _n( k) \Big( \frac{n^2}{k^2}  \sum  _{i \in
I}\eg_i(\ebo)\eg^t_i(\ebo)  +\frac{n^2 }{(n-k)^2} 
\sum  _{j \in J} \eg_j(\ebo)\eg^t_j(\ebo)
\Big) \hat{\el} _n( k)\Big]
 \nonumber \\  && 
 + \Big[  2\hat{\el}^t _n( k)
  \Big((\em_n^0)^{ \frac{1}{2}} (\textbf{V}^0_{1n})^{-1}\frac{n}{k}
\sum  _{i \in I}\egg_i(\ebo) - (\em_n^0)^{ \frac{1}{2}} (\textbf{V}^0_{2n})^{-1}\frac{n}{n-k}  \sum  _{j \in J} \egg_j(\ebo) \Big )  (\hat{\eb}_n( k)-\ebo) \Big ] 
\nonumber \\  &&  
 - \Big[ 2\hat{\el}^t_n( k) \Big(
   \frac{n}{k}  \sum  _{i \in I}\eg_i( \ebo)\frac{\partial ( \em_n^{ \frac{1}{2}}(\eb) (\eV_{1n}(\eb))^{-1})}{\partial \eb}
   %\nonumber \\  &&  
   + \frac{n}{n-k}  \sum  _{j \in J}\eg_j( \ebo)\frac{\partial ( \em_n^{ \frac{1}{2}}(\eb) (\eV_{2n}(\eb))^{-1})}{\partial \eb}
   \Big ) (\hat{\eb}_n( k)-\ebo) \Big] 
  \nonumber \\  && 
  + \frac{1}{3 !} \Big [T_1+ 3 T_2 + 3 T_3 +T_4    \Big], 
\end{eqnarray}
where \\
$ T_1 = \sum  _{r=1}^d  \sum _{l=1}^d \sum  _{s=1}^d \frac{\partial^3 {\cal Z}_{nk}(\del_0,\el^{(1)}_{rsl},\eb^{(1)}_{rsl})}{\partial \beta_r \partial \beta_s \partial \beta_l } (\hat\beta_{n,r}-\beta^0_r) (\hat\beta_{n,s}-\beta^0_s)(\hat\beta_{n,l}-\beta^0_l)$,  \\
$ T_2 = \sum  _{r=1}^d  \sum  _{l=1}^d \sum  _{s=1}^d \frac{\partial^3 {\cal Z}_{nk}(\del_0,\el^{(2)}_{rsl},\eb^{(2)}_{rsl})}{\partial \lambda_r \partial \lambda_s \partial\beta_l } (\hat\lambda_{n,r}) (\hat\lambda_{n,s})(\hat\beta_{n,l}-\beta^0_l)$,\\
$ T_3 =  \sum _{r=1}^d  \sum  _{l=1}^d \sum  _{s=1}^d \frac{\partial^3 {\cal Z}_{nk}(\del_0,\el^{(3)}_{rsl},\eb^{(3)}_{rsl})}{\partial \lambda_r \partial \beta_s \partial \beta_l } (\hat\lambda_{n,r}) (\hat\beta_{n,s}-\beta^0_s)(\hat\beta_{n,l}-\beta^0_l)$,\\
$T_4 = \sum  _{r=1}^d  \sum  _{l=1}^d \sum  _{s=1}^d \frac{\partial^3 {\cal {Z}}_{nk}(\del_0,\el^{(4)}_{rsl},\eb^{(4)}_{rsl})}{\partial \lambda_r \partial \lambda_s \partial \lambda_l } (\hat\lambda_{n,r}) (\hat\lambda_{n,s})(\hat\lambda_{n,l})$,  \\ \\
%%%%%%%%%%%%%%%%%%%%%%%%%%%%%%%%%%%%%%%%%%%%%%%%%%%%%%%%%%%%%%ù
for $1 \leq r \leq d$, $\hat \beta_{n,r}$ being the r-th component of $\hat{\eb}_n(k)$, $\hat \lambda_{n,r}$ is the r-th component of $\hat{\el}_n(k)$. For all $1 \leq r,s,l \leq d$, $\el^{(b)}_{rsl}= u^{(b)}_{rsl} (\hat{\eb}_n(k)-\ebo)$ and $\eb^{(b)}_{rsl}= \ebo+ v^{(b)}_{rsl} (\hat{\eb}_n(k)-\ebo)$, with $u^{(b)}_{rsl},v^{(b)}_{rsl} \in [0,1]$ and $b \in \{1,2,3,4 \}$. \\ \\
We note that, the derivatives $\partial ( \em_n^{ \frac{1}{2}}(\eb) (\eV_{1n}(\eb))^{-1}) ) / \partial \eb$ and $\partial ( \em_n^{ \frac{1}{2}}(\eb) (\eV_{2n}(\eb))^{-1}) ) / \partial \eb$ are considered term by term. \\
Now, we replace $\hat{\el}_n(k)$ in the relation (\ref{DLz}) by their value obtained in $(\ref{el1})$. For the first term of (\ref{DLz}), using notations given above, we find that this term is equal to $2(n\epsi ^0_n )^t (\es _n^0+\eephi _n ^0)^{-1}\epsi^0_n$. \\
Similarly, the second term of (\ref{DLz}) is $n(\epsi_n^0)^t (\es ^0_n+\eephi _n ^0)^{-1}\epsi^0_n$. We know that, $\eV_{1n}^0=k^{-1}
\sum  _{i \in I}\egg_i(\ebo) $ and $\eV_{2n}^0=(n-k)^{-1}
\sum  _{j \in J}\egg_j(\ebo) $. Then the third term of (\ref{DLz}) converge almost surely to zero, as $n\rightarrow \infty$. \\
By the central limit theorem, we have that $k^{-1}\sum  _{i \in I}\eg_i(\ebo)$ = $O_{\eP}(k^{-1/2})$ and  $(n-k)^{-1} \sum  _{j \in J}\eg_j(\ebo)$ = $O_{\eP}((n-k)^{-1/2})$. Then, the fourth term of (\ref{DLz}) is $o_{\eP}(n (\epsi_n^0)^t (\es_n ^0+\eephi_n ^0)^{-1}\epsi_n^0)$. \\
For the last term of (\ref{DLz}), using assumptions (A2)-(A4) and by elementary calculations, we prove that this term is $\xi_n (k)$, where 
\begin{equation}
\label{xi}
\xi_n (k) \equiv o_{\eP}(\|\hat\eb_n( k) -\ebo\|_2)+o_{\eP}(\|\hat \el_n ( k)\|_2)+o_{\eP}(\|\hat \el_n (k)\|_2 \|\hat \eb_n ( k)-\ebo\|_2).
\end{equation}
Combining the obtained results, we obtain
\begin{equation}
\label{ez1}
{\cal Z}_{nk}(\del_0,\hat \el_n ( k),\hat  \eb_n (k))=n(\epsi^0 _n)^t (\es^0 _n+\eephi^0 _n)^{-1}\epsi^0 _n(1+o_{\eP}(1))+ \xi_n (k).
\end{equation}
We can see that, when $ \eb-\ebo=\textbf{0}_d$, i.e, $ \eb=\ebo$, $n(\epsi^0 _n)^t (\es^0 _n+\eephi^0 _n)^{-1}\epsi^0 _n$ achieves its maximum in the neighborhood of $\ebo$. This means when $n\rightarrow \infty$, ${\cal Z}_{nk}(\el,\del_0,\eb)$ has one local maximum in any $\eta$-neighborhood of $\ebo$ and then $\hat{\eb}_n(k)- \eb_0=o_{\eP}(1)$. \\
For $\es^0 _n$, by the law of large numbers, the terms $(n/k^2) \sum   _{i\in I}\et_i(\ebo) \et^t_i(\ebo)$ and $(n/(n-k)^2) \sum   _{j\in J}\et_j(\ebo)\et^t_j(\ebo)$ converge almost surely to the expected of $\et_i(\ebo)\et^t_i(\ebo)$ and $\et_j(\ebo)\et^t_j(\ebo)$, respectively as $n\rightarrow \infty$. Using this fact and the fact that $-\eV_{1n}^0 \overset{a.s}{\longrightarrow}\eV$ and $-\eV_{2n}^0 \overset{a.s}{\longrightarrow} \eV$, by the law of large numbers, we can demonstrate easily that $\es^0 _n=\textbf{I}_d+O_{\eP}(n^{-1/2})$. Using this fact and relation (\ref{el1}), we obtain that $\hat{\el}_n (k)= \epsi ^0 _n+ o_{\eP}(n^{-1/2})$. By the central limit theorem and the fact that $\eE[g_i(\ebo)]=0$ for $i=1,\ldots,n$, each term of $\epsi^0 _n$ is $O_{\eP}(n^{-1/2})$, which implies $\hat{\el}_n(k)=O_{\eP}(n^{-1/2})$.
The lemma is completely proved.\hspace*{\fill}$\blacksquare$  \ \\ \ \\
%\label{prrof theo1}
\noindent {\bf \textbf{ Proof of Theorem \ref{theo1}}}. Using the proof of Proposition \ref{prop1}, we have $\hat{\el}_n( k)= \epsi^0_n + o_{\eP}(n^{-1/2})$ and $\es^0 _n=\textbf{I}_d+O_{\eP}(n^{-1/2})$. By the Linderberg-Feller Theorem, we have that $\sqrt n \hat{\el}_n ( k)  \overset{{\cal L}} {\underset{n \rightarrow \infty}{\longrightarrow}} {\cal N}(0,\textbf{I}_d)$. Then, using also  relation (\ref{el1}), we obtain
\begin{eqnarray*}
\label{ez11}
{\cal Z}_{nk}(\del_0,\hat \el_n (k),\hat \eb_n ( k))&=&n(\epsi^0_n)^t (\es^0_n)^{-1}\epsi^0_n +o_{\eP}(1)+ \xi_n (k)\nonumber \\  && = n(\epsi^0_n)^t [I_d+O_{\eP}(n^{-1/2})]^{-1} \epsi^0_n+o_{\eP}(1) + \xi_n (k)\nonumber \\  && = n(\epsi^0_n)^t \epsi^0_n +o_{\eP}(1) +\xi_n (k)\nonumber \\  && = n\hat{\el}^t_n( k) \hat{\el}_n (k)+o_{\eP}(1)  +\xi_n (k)\nonumber \\  && = \sqrt n \hat{\el} ^t_n(k) \hat{\el}_n (k)\sqrt n+o_{\eP}(1)+\xi_n (k),
\end{eqnarray*}
where $\xi_n (k)$ is given by the relation (\ref{xi}). Then, since $\sqrt n \hat{\el}_n ( k) \sim {\cal N}(0,\textbf{I}_d)$, we deduce
\begin{equation*}
{\cal Z}_{nk}(\del_0,\hat \el_n (k),\hat \eb_n (k))\overset{{\cal L}} {\underset{n \rightarrow \infty}{\longrightarrow}} \chi^2(d).
\end{equation*}
The theorem is proved.\hspace*{\fill}$\blacksquare$  \ \\ \ \\
\noindent {\it \textbf{ Proof of Proposition \ref{prop2}}}. For each method, the proof is similar to the proof of Proposition \ref{prop1}. \hspace*{\fill}$\blacksquare$  \ \\ \ \\
\noindent {\it \textbf{ Proof of Theorem \ref{theo2}}}. Using Proposition \ref{prop2}, the proof of this theorem is similar to the proof of Theorem \ref{theo1}, for the three proposed methods. \hspace*{\fill}$\blacksquare$ 

\section*{References}
%%%%%%%%%%%%%%%%%%%%%%%
Bai, J., 1999. Likelihood ratio tests for multiple structural changes. \textit{Journal of Econometrics}. 91:299

-323. \\
%%%%%%%%%%%%%%%%%%%%%%%
Boldea, O., Hall, A.R., 2013. Estimation and inference in unstable nonlinear least square models.

\textit{Journal of}\textit{Econometrics}, {172}(1), 158-167.\\
%%%%%%%%%%%%%%%%%%%%%%%
Chen, S., 1993. On the accuracy of empirical likelihood confidence region for linear regression model.

 \textit{Ann. Inst. Statistic. Math}, {45}, 621-637.\\
%%%%%%%%%%%%%%%%%%%%%%%
Chen, S., 1994. Empirical likelihood confidence interval for linear regression coefficients. \textit{Ann.Inst. Stati-}

\textit{-stic. Math}, {45}, 621-637.\\
 %%%%%%%%%%%%%%%%%
Chen, S., Cui, H., 2003. An extended empirical likelihood for generalized linear models.
\textit{Statistica Sinicia }

 {13}, 69-81.\\
 %%%%%%%%%%%%%%%%%
%Ciuperca, G., 2011. {A general criterion to determinate the number of change-points}. \textit{Statistics and}
%
% \textit{Probability Letters}, {81}(8), 1267-1275.\\
 %%%%%%%%%%%%%%%%%
Ciuperca G., 2011. Empirical likelihood for nonlinear model with missing responses.
\textit{Journal of Statist-}

\textit{-ical Computation and Simulation}, {83}(4), 737-756.
\\
 %%%%%%%%%%%%%%%%%
Ciuperca, G., Salloum, Z., 2013. {Empirical likelihood test in a posteriori change-point
nonlinear model}. 

\textit{Accepted to Metrika}.\\
 %%%%%%%%%%%%%%%%%
  %%%%%%%%%%%%%%%%%
Horvitz, DG., Thompson, DJ., 1952. A generalization of sampling without replacement from a finite

 universe. {\textit{J Am Stat assoc}},  {47}, 663-685.\\
 %%%%%%%%%%%%%%%%%
Jing, B., 1995. {Two-sample empirical likelihood method}. \textit{Statistics and Probability Letters},
 {24}, 315-319.\\
 %%%%%%%%%%%%%%%%%
 %%%%%%%%%%%%%%%%%
Kim, H.j., Siegmund, D., 1989.
{The likelihood ratio test for a change-point in simple linear regression 

model}. \textit{Biometrika},  {76}, 409-423.\\
%%%%%%%%%%%%%%%%%
%%%%%%%%%%%%%%%%%
%%%%%%%%%%%%%%%%%
Kolaczyk, E. D., 1994. {Empirical likelihood for generalized linear models}. {\textit{Statistica Sinicia}}, {4}, 199-218.\\
%%%%%%%%%%%%%%%%%
Little, R.J.A, Rubin, D.B., 1987. {Statistical analysis with missing data}. \textit{John Wiley and sons, Inc., New}

 \textit{York}.\\
  %%%%%%%%%%%%%%%%%
Liu, Y., Zou, C., Zhang, R.,2008. {Empirical likelihood for the two-sample mean problem}. \textit{Statistics and }

\textit{Probability Letters},
 {78}, 548-556.\\
%%%%%%%%%%%%%%%%%
Liu, Y.,  Zou,  C.,  Zhang, R., 2008.
Empirical likelihood ratio test for a change-point in linear regression 

model. \textit{Communications in Statistics-Theory and Methods},  {37}, 2551-2563.\\
%%%%%%%%%%%%%%%%%
%%%%%%%%%%%%%%%%%
Muller, U.U., 2009. Estimating linear functionals in nonlinear regression with responses missing at ra-

-ndom. \textit{Annal statistic},  {37}(5A), 2245-2277.\\
%%%%%%%%%%%%%%%%%
%%%%%%%%%%%%%%%%%
Ning, W., Pailden, J., Gupta, A., 2012.
 {Empirical likelihood ratio test for the epidemic change model}. 
 
 \textit{Journal of Data Science},  {10}, 107-127.\\
%%%%%%%%%%%%%%%%%
%%%%%%%%%%%%%%%%%
Owen, A. B., 1988. {Empirical likelihood ratio confidence intervals for a single functional}.  {\textit{Biometrika}}, 

{75}, 237-249.\\
%%%%%%%%%%%%%%%%%
%%%%%%%%%%%%%%%%%
Owen, A. B., 1990. {Empirical likelihood ratio confidence regions}. {\textit{Annal statistic}}, {18}, 90-120.\\
%%%%%%%%%%%%%%%%%
Owen, A. B., 1991. {Empirical likelihood for linear model}. {\textit{Annal statistic}}, {19}, 1725-1747.\\
%%%%%%%%%%%%%%%%%
Qin, Y., Li, L., Lei, Q., 2009. {Empirical likelihood for linear regression model with missing responses}. 

 {\textit{Statistics and Probability Letters}}, {79}, 1391-1396.\\
%%%%%%%%%%%%%%%%%
Seber G., Wild C., 2003.
{Testing for structural change in regression quantiles, Wiley series in probability 

and mathematical Statistics, Wiley, Hoboken, NJ}.\\
%%%%%%%%%%%%%%%%%%%%%%
Sun, Z., Wang, Q., Dai, P., 2009. model checking for partially linear models with missing responses at

 random. \textit{Journal of Multivariate Analysis}, { 100}, 636-651.\\
%%%%%%%%%%%%%%%%%%%%%%%%%%%%%
%%%%%%%%%%%%%%%%%
Wei, Yu., Cuizhen, Ni., Wangli, Xu., 2013. {An empirical likelihood inference for the coefficient difference 

of a two-sample linear model with missing response data}. \textit{Metrika},
 {81}(8), 1267-1275.\\
 %%%%%%%%%%%%%%%%%
%%%%%%%%%%%%%%%%%
Xue, L., 2009. {Empirical likelihood for linear models with missing responses}.  \textit{Journal of Multivariate}

 \textit{Analysis}, { 100}, 1353-1366.\\
%%%%%%%%%%%%%%%%
%%%%%%%%%%%%%%%%%
Zi, X., Zou, C., Liu, Y., 2010.
{Two-sample empirical likelihood method for difference between coefficients 

in linear regression model}.
 \textit{Stat Papers}.\\
 %%%%%%%%%%%%%%%%%
%%%%%%%%%%%%%%%%%
%%%%%%%%%%%%%%%%%
Zou C., Liu Y., Qin P., Wang Z., 2007. {Empirical likelihood ratio test for a change point}.\textit{Statistics and }

\textit{Probability Letters}, { 77}, 374-382.\\
%%%%%%%%%%%%%%%%%

\end{document}